\documentclass{amsart}
\usepackage{graphicx}
\usepackage{xspace,colortbl}
\usepackage{amssymb}
\theoremstyle{definition}

\theoremstyle{remark}

\numberwithin{equation}{section}

\begin{document}

\title[Some estimates for the bilinear fractional integrals]{Some estimates for the bilinear fractional integrals on the Morrey space}

\author[Xiao Yu \ \ Xiangxing Tao \ \  Huihui Zhang and Jianmiao Ruan]{{Xiao Yu \ \ Xiangxing Tao   \ \ Huihui Zhang  and Jianmiao Ruan}}%

\thanks{Supported by the National Natural Science Foundation of China\,(11561057, 11171306),
  the Natural Science Foundation of Jiangxi Province\,(20151BAB211002)  and  the Science Foundation of Jiangxi Education Department\,(GJJ151054). }%
\thanks{{\bf Xiao Yu: } Department of Mathematics, Shangrao Normal University, Shangrao 334001, P. R. China. E-mail: yx2000s@163.com; \ \ xyuzju@163.com}
\thanks{{\bf Xiangxing Tao: } Department of Mathematics, Zhejiang University of Science and Technology, Hangzhou 310023,  P.R.China. E-mail:xxtao@zust.edu.cn}
\thanks{{\bf Huihui Zhang: } Department of Mathematics, Shangrao Normal University, Shangrao 334001, P. R. China. E-mail:zhanghuihuinb@163.com}
\thanks{{\bf Jianmiao Ruan: } Department of Mathematics, Zhejiang International Studies University, Hangzhou 310012, P. R. China. E-mail:rjmath@163.com}


\date{August 10, 2018}
\begin{abstract}
In this paper, we are interested in the following bilinear fractional integral operator $B\mathcal{I}_\alpha$ defined by

\[
B\mathcal{I}_{\alpha}({f,g})(x)=\int_{%
\mathbb{R}
^{n}}\frac{f(x-y)g(x+y)}{|y|^{n-\alpha}}dy,
\]
with $0< \alpha<n$. We prove the weighted  boundedness of $B\mathcal{I}_\alpha$ on the  Morrey type spaces. Moreover, an Olsen type inequality for $B\mathcal{I}_\alpha$ is also given.

 \end{abstract}
\maketitle

\section{\textbf{Introduction}}
In 1992, Grafakos \cite{Gra} studied the multilinear fractional integral operator $\mathcal{I}_{\alpha ,\vec{\theta }}$ with its definition defined by
\[
\mathcal{I}_{\alpha ,\vec{\theta }}(\vec{f})(x)=\int_{%
\mathbb{R}
^{n}}\frac{1}{\left\vert y\right\vert ^{n-\alpha }}\prod%
\limits_{i=1}^{m}f_{i}(x-\theta _{i}y)dy,
\]%
where \
\[
\vec{f}=(f_1, \cdots, f_m)
\]%
and \
\[
\vec{\theta} =\mathbf{(}\theta _{1},\theta _{2},...,\theta _{m})\text{ \ }
\]%
is a fixed vector with distinct nonzero components.

For a special case of $\mathcal{I}_{\alpha ,\vec{\theta }}$,  the following bilinear fractional integral was also studied by  Kenig and Stein in  \cite{KS}.
\[
B\mathcal{I}_{\alpha}({f,g})(x)=\int_{%
\mathbb{R}
^{n}}\frac{f(x-y)g(x+y)}{|y|^{n-\alpha}}dy, \ \ \ 0<\alpha<n.
\]

As the operator $B\mathcal{I}_{\alpha}$ is similar to the bilinear Hilbert transform if we take $\alpha\rightarrow 0$, many authors pay much attention to such operator and they proved the  boundedness of $B\mathcal{I}_\alpha$ on the product function spaces.  One may see \cite{BDMT,BMMN,ChF1,ChF2,ChF3,DL,WC,ZC} et al. for more details.

Meanwhile, it is well known that in the last 70s,  Muckenhoupt and Wheeden (\cite{M1, MW}) introduced the $A_p$ weight class  and $A_{(p,q)}$ weight class which are very adopted for the weighted estimates of the singular integrals and fractional integrals. Now, let us introduce the  definitions of $\omega\in A_p$ and $\omega\in A_{(p,q)}$ respectively.

\textbf{Definition 1.1.} (\cite{M1}) We say a non-negative function $\omega(x)$ belongs to the Muckenhoupt class $A_p$ with $1<p<\infty$ if
$$[\omega]_{A_p}:=\sup\limits_Q\left(\frac{1}{|Q|}\int_Q\omega(x)dx\right)\left(\frac{1}{|Q|}\int_Q\omega(x)^{1-p'}dx\right)^{p-1}< \infty \eqno(1.1)$$
for any cube $Q$ and $1/p+1/p'=1$.

In case $p=1$, $\omega\in A_1$ is understood as there exists a positive constant $C$ such that
$$\frac{1}{|Q|}\int_Q\omega(y)dy\leq C\omega(x) \eqno(1.2)$$
for a.e. $x\in Q$ and any cube $Q$. For the case $p=\infty$,  we define  $A_{\infty}=\bigcup\limits_{1<p<\infty}A_p$.

\textbf{Definition 1.2.} (\cite{MW}) We say that a non-negative function $\omega(x)$  belongs to $A_{(p,q)}$ weight class with $1<p<q<\infty$ if
$$[\omega]_{A_{p,q}}:=\sup\limits_Q\left(\frac{1}{|Q|}\int_Q\omega(x)^qdx\right)^{1/q}\left(\frac{1}{|Q|}\int_Q\omega(x)^{-p'}dx\right)^{1/p'}<\infty.\eqno(1.3)$$

Since the last 90s, the multilinear  theory for the singular integral operators was developed a lot. For example,  in 2002, Grafakos and Torres \cite{GT} introduced the multilinear C-Z theory. Later,  Lerner et al. \cite{LOPTT} introduced a new kind of multiple weight which is very adopted for the weighted norm inequalities of the multilinear singular integrals.  Following their work, Chen and Xue \cite{CX}, as well as Moen independently \cite{Mo1}, introduced a new type of multiple fractional type $A_{(\vec{P},q)}$ weight class .   Now, let us give the definition of $A_{(\vec{P},q)}$ weight class.

\textbf{Definition 1.3.} (\cite{CX, Mo1}) Let $1\leq
p_1,\cdots,p_m, 1/p=1/p_1+\cdots+1/p_m$ and $q>0$. Suppose that
$\vec{\omega}=(\omega_1,\cdots,\omega_m)$ and each $\omega_i$ is a
nonnegative function on $\Bbb{R}^n$. We say that $\vec{\omega}\in A_{(\vec{p},q)}$ if
it satisfies
$$
\begin{array}{ll}
[\vec{\omega}]_{A_{(\vec{P}, q)}}:=\sup\limits_{Q}\left(\frac{1}{|Q|}\int_Q\nu_{\vec{\omega}}^q(x)dx\right)^{1/q}\prod\limits_{i=1}^m\left(\frac{1}{|Q|}\int_Q\omega_i^{-p_i'}(x)dx\right)^{1/p_i'}<\infty,
\end{array}
 \eqno (1.4)
$$
where $\nu_{\vec{\omega}}=\prod\limits_{i=1}^m\omega_i$.  Moreover, for the case  $p_i=1$,
$\left(\frac{1}{|Q|}\int_Q\omega_i^{-p_i'}\right)^{1/p_i'}$ is
understood as $(\inf\limits_Q\omega_i)^{-1}$.\\

Chen and Xue, as well as Moen independently, proved the following theorem.

\textbf{Theorem A.} (\cite{CX,Mo1}) Suppose that $0<\alpha<mn$, $1<p_1,\cdots,p_m<\infty$. If $1/p=\sum\limits_{i=1}^m1/p_i$ and $1/q=1/p-\alpha/n$. Then, $\vec{\omega}\in A_{(\vec{P},q)}$ if and only if  the following multiple weighted norm inequalities holds:

$$\|\mathcal{I}_{\alpha,m}(\vec{f})\|_{L^q(\nu_{\vec{\omega}}^q)}\leq C\prod\limits_{i=1}^m\|f_i\|_{L^{p_i}(\omega_i^{p_i})}.$$
Here, $\mathcal{I}_{\alpha,m}$ denotes the multilinear fractional integral operator and its definition can be stated as
$$\mathcal{I}_{\alpha,m}(\vec{f})(x)=\displaystyle\int_{(\Bbb{R}^n)^m}\frac{f_1(y_1)f_2(y_2),\cdots,f_m(y_m)}{(|x-y_1|+\cdots+|x-y_m|)^{mn-\alpha}}dy_1dy_2\cdots
dy_m.$$

For the study of the weighted theory for $B\mathcal{I}_\alpha$ with the multiple fractional type weight class,  Hoang and Moen \cite{HM,M} did some excellent work to show that the operator $B\mathcal{I}_\alpha$ satisfy several weighted estimates on the product $L^p$ spaces.

On the other hand, in order to
study the local behavior of solutions to second order elliptical partial
differential equations, Morrey \cite{Morrey} introduced the Morrey space. The Morrey space \ $\mathcal{M}_{q}^{p}(%
\mathbb{R}
^{n}),$ $0<q\leq p<\infty ,$ \ is the collection of all measurable functions
\ $f$ \ with its definition defined by
\[
\mathcal{M}_{q}^{p}(\Bbb{R}^n):=\{f\in \mathcal{M}_{q}^{p}(\Bbb{R}^n):\left\Vert f\right\Vert _{\mathcal{M}_{q}^{p}(%
\mathbb{R}
^{n})}=\sup_{\substack{Q\subset \Bbb{R}^n\\Q:\text{cubes}}}\left\vert Q\right\vert ^{1/p-1/q}\left\Vert f\chi
_{Q}\right\Vert _{L^{q}(%
\mathbb{R}
^{n})}<\infty\}.
\]%
Many authors  studied the weighted norm inequalities for integral operators on the Morrey type spaces, readers may see \cite{I1,I2,ISST1,ISST2,KoS} et al. or the summary article \cite{I3} to find more details. Here we would like to mention that in \cite{I1,I2,ISST1}, Iida et al. introduced  the following  new fractional type multiple weight condition as follows.

$$
\begin{array}{ll}
[\vec{\omega}]_{q_0,q,\vec{P}}&:=\sup\limits_{\substack{Q\subset Q'\\Q,Q':\text{cubes}}}\left(\frac{|Q|}{|Q'|}\right)^{1/q_0}\left(\frac{1}{|Q|}\int_Q(\omega_1(x)\omega_2(x))^qdx\right)^{1/q}\\
&\times\prod\limits_{i=1}^2\left(\frac{1}{|Q|}\int_{Q'}\omega_i(y_i)^{-p_i'}dy_i\right)^{1/p_j'}<\infty.
\end{array}
\eqno(1.5)
$$

Iida et al. \cite{I1,I2,ISST1} proved  that the above multiple weight condition is very adopted for the weighted norm inequalities of the operator $\mathcal{I}_\alpha$ on the   Morrey type  space and they proved the following theorem.

\textbf{Theorem B.} (\cite{I1,ISST1}) Let $0<\alpha<mn$, $1<p_1,\cdots,p_m<\infty$, $1/p=\sum\limits_{i=1}^m1/p_i$.  Then,we assume that $0<p\leq p_0<\infty$ and $0<q\leq q_0<\infty$ with $1/q_0=1/p_0-\alpha/n$ and $q/q_0=p/p_0$. Moreover, for $\vec{f}=(f_1,\cdots,f_m)$ and $\vec{\omega}=(\omega_1,\omega_2)$, we denote
$$\|\vec{f}\|_{\mathcal{M}^{p_0}_{\vec{P}}}:=\sup\limits_{Q\subset \Bbb{R}^n}|Q_0|^{\frac{1}{p_0}}\prod\limits_{i=1}^m\left(\frac{1}{|Q|}\int_Q|f_i(y_i)|^{p_i}dy_i\right)^{\frac{1}{p_i}},$$
and
$$\nu_{\vec{\omega}}(x)=\prod\limits_{i=1}^m\omega_i(x).$$ If there exist $a>1$ satisfying
 $$[\vec{\omega}]_{aq_0,q,\vec{P}}<\infty$$
where $\vec{P}=(p_1,\cdots,p_m)$ and $a>1$, then  there exist a positive constant $C$ independent of $f_i$, such that
$$\|\mathcal{I}_{\alpha,m}(\vec{f})\nu_{\vec{\omega}}\|_{\mathcal{M}^{q_0}_q}\leq C\|(f_1\omega_1,\cdots,f_m\omega_m)\|_{\mathcal{M}^{p_0}_{\vec{P}}}.$$

In \cite{HY}, He and Yan studied the weighted bounedness of  $B\mathcal{I}_\alpha$ on $\mathcal{M}_{q}^{p}(\Bbb{R}^n)$ with $0<q<1$. Thus,  it is natural to ask whether we can prove the weighted norm inequalities for $B\mathcal{I}_\alpha$ on $\mathcal{M}_{q}^{p}(\Bbb{R}^n)$ with $q>1$? In this paper, we will give a positive answer to this question.

Motivated by the above backgrounds, in this paper, we will give the weighted boundedness of $B\mathcal{I}_\alpha$ on the Morrey type space with the fractional type multiple weights condition proposed by Iida et al.  Our results can be stated as follows.

\textbf{Theorem 1.1.} Suppose $0<\alpha<n$, $p_1>r>1$, $p_2>s>1$, $1/r+1/s=1$, $1/p=1/p_1+1/p_2$, $1<p_1,p_2<\infty$,  $0<p\leq p_0<\infty$, $0<q\leq q_0<\infty$. Let
$$1/q_0=1/p_0-\alpha/n, \ \ q/q_0=p/p_0 \ \ \text{and}\ \  \nu_{\vec{\omega}}(x)=\prod\limits_{i=1}^2\omega_i(x).$$
Moreover, assume that either $p$ or $q$ satisfies one of the following condition:
$$p>1\ \ \ \ \ \text{or} \ \ \ \ \ q>\frac{1}{2}.$$
If there exists $a>1$, such that
$[\vec{\omega}]_{aq_0,q,(\frac{sp_1}{s+p_1},\frac{rp_2}{r+p_2})}<\infty$, that is
$$
\begin{array}{ll}
\sup\limits_{\substack{Q\subset Q'\\Q,Q':\text{cubes}}}&\left(\frac{|Q|}{|Q'|}\right)^{\frac{1}{aq_0}}\left(\frac{1}{|Q|}\int_Q\nu_{\vec{\omega}}(x)^qdx\right)^{1/q}\left(\frac{1}{|Q'|}\int_{Q'}\omega_1(x)^{-\frac{p_1r}{p_1-r}}\right)^{1/r-1/p_1}\\
&\times\left(\frac{1}{|Q'|}\int_{Q'}\omega_2(x)^{-\frac{p_2s}{p_2-s}}\right)^{1/s-1/p_2}<\infty.
\end{array}
\eqno(1.6)
$$
Then, there exists a positive constant $C$ independent of $f$ and $g$, such that
$$\|B\mathcal{I}_{\alpha}(f,g)\nu_{\vec{\omega}}\|_{\mathcal{M}^{q_0}_{q}}\leq C[\vec{\omega}]_{aq_0,q,(\frac{sp_1}{s+p_1},\frac{rp_2}{r+p_2})}\|(f\omega_1,g\omega_2)\|_{\mathcal{M}^{p_0}_{\vec{P}}}.\eqno(1.7)$$
\medskip

\textbf{Remark 1.2.} Note that for the operator $\mathcal{I}_{\alpha,2}$
$$\mathcal{I}_{\alpha,2}(\vec{f})(x)=\displaystyle\int_{(\Bbb{R}^n)^m}\frac{f_1(y_1)f_2(y_2)}{(|x-y_1|+|x-y_2|)^{2n-\alpha}}dy_1dy_2.$$
As mentioned in \cite[p.629]{M}, if we denote $\delta$   is the point mass measure at the origin, then  we know that the kernel of $\mathcal{I}_{\alpha,2}$,
$$K_\alpha(u,v)=(|u|+|v|)^{-2n+\alpha}$$
has a singularity at the origin in $\Bbb{R}^{2n}$ as opposed to the kernel of  $B\mathcal{I}_{\alpha}$
$$k_{\alpha}(u,v)=\frac{\delta(u+v)}{|u|^{n-\alpha}},$$
which has a singularity along a line.  Thus, we   conclude that  Theorem 1.1 parallel earlier results by the authors \cite{ISST1} for the  less singular bilinear fractional integral operator $\mathcal{I}_{\alpha,2}$.
\medskip

\textbf{Remark 1.3.} If we take $p=p_0$ and $q=q_0$, we obtain the following result  proved by Hoang and Moen \cite{HM}.
\medskip

\textbf{Corollary 1.4.} (\cite{HM}) Suppose that there exist real numbers $\alpha ,p_1, r , p_2, s, p$ and $q$ satisfying the same conditions as in Theorem 1.1. If $1/p_1+1/p_2-1/q=\alpha/n$ and $\vec{\omega}\in A_{((\frac{p_1s}{p_1+s},\frac{p_2r}{p_2+r}),q)}$, then there exists a positive constant $C$ independent of $f$ and $g$, such that
$$\|B\mathcal{I}_\alpha(f,g)\|_{L^q(\nu_{\vec{\omega}}^q)}\leq C\|f\|_{L^{p_1}(\omega_1^{p_1})}\|g\|_{L^{p_2}(\omega_2^{p_2})}.\eqno(1.8)$$

\textbf{Proof of Corollary 1.4.} By the definition of the Morrey space, it suffices to show that $$[\vec{\omega}]_{aq_0,q,(\frac{sp_1}{s+p_1},\frac{rp_2}{r+p_2})}<\infty\ \  (q=q_0).\eqno(1.9)$$

In fact, as $\vec{\omega}\in A_{((\frac{p_1s}{p_1+s},\frac{p_2r}{p_2+r}),q)}$, we have $\nu_{\vec{\omega}}^q=\prod\limits_{i=1}^2\omega_i^q\in A_{2q}$ or $\nu_{\vec{\omega}}^q\in A_{1+q(1-1/p)}$. Then, we know that $\nu_{\vec{\omega}}^q$ satisfies the reversed H\"{o}lder inequality (see Section 2).  That is, if we choose $a=1+\epsilon$ where $\epsilon \in \Bbb{R}^+$ and $\epsilon$ is small enough, there is
$$\left(\frac{1}{|Q|}\int_Q\nu_{\vec{\omega}}(x)^{aq}dx\right)^{1/aq}\leq C\left(\frac{1}{|Q|}\int_Q\nu_{\vec{\omega}}(x)^qdx\right)^{1/q}.$$
Recall that $q=q_0$, we may have
\begin{align*}
&\left(\frac{|Q|}{|Q'|}\right)^{\frac{1}{aq_0}}\left(\frac{1}{|Q|}\int_Q(\omega_1(x)\omega_2(x))^qdx\right)^{1/q}\left(\frac{1}{|Q'|}\int_{Q'}\omega_1(x)^{-\frac{p_1r}{p_1-r}}dx\right)^{1/r-1/p_1}\\
&\times\left(\frac{1}{|Q'|}\int_{Q'}\omega_2(x)^{-\frac{p_2s}{p_2-s}}dx\right)^{1/s-1/p_2}\\
&\leq \left(\frac{|Q|}{|Q'|}\right)^{\frac{1}{aq_0}}\left(\frac{1}{|Q|}\int_Q(\omega_1(x)\omega_2(x))^{aq}dx\right)^{1/aq}\left(\frac{1}{|Q'|}\int_{Q'}\omega_1(x)^{-\frac{p_1r}{p_1-r}}dx\right)^{1/r-1/p_1}\\
&\times\left(\frac{1}{|Q'|}\int_{Q'}\omega_2(x)^{-\frac{p_2s}{p_2-s}}dx\right)^{1/s-1/p_2}\\
&=\left(\frac{1}{|Q'|}\int_Q(\omega_1(x)\omega_2(x))^{aq}dx\right)^{1/aq}\left(\frac{1}{|Q'|}\int_{Q'}\omega_1(x)^{-\frac{p_1r}{p_1-r}}dx\right)^{1/r-1/p_1}\\
&\times\left(\frac{1}{|Q'|}\int_{Q'}\omega_2(x)^{-\frac{p_2s}{p_2-s}}dx\right)^{1/s-1/p_2}\\
&\leq\left(\frac{1}{|Q'|}\int_{Q'}(\omega_1(x)\omega_2(x))^{aq}dx\right)^{1/aq}\left(\frac{1}{|Q'|}\int_{Q'}\omega_1(x)^{-\frac{p_1r}{p_1-r}}dx\right)^{1/r-1/p_1}\\
&\times\left(\frac{1}{|Q'|}\int_{Q'}\omega_2(x)^{-\frac{p_2s}{p_2-s}}dx\right)^{1/s-1/p_2}\\
&\leq\left(\frac{1}{|Q'|}\int_{Q'}(\omega_1(x)\omega_2(x))^qdx\right)^{1/q}\left(\frac{1}{|Q'|}\int_{Q'}\omega_1(x)^{-\frac{p_1r}{p_1-r}}dx\right)^{1/r-1/p_1}\\
&\times\left(\frac{1}{|Q'|}\int_{Q'}\omega_2(x)^{-\frac{p_2s}{p_2-s}}dx\right)^{1/s-1/p_2}<\infty,
\end{align*}
where  the last inequality follows from the reversed H\"{o}lder inequality for $\nu_{\vec{\omega}}^q=(\omega_1\omega_2)^q$ and we obtain (1.9).

\textbf{Remark 1.5.} For the case $0<q<1$ in Theorem 1.1,  our result is  still different from \cite[Theorem 4.6]{HY}.

\section{Preliminaries}
In this section, we will give some lemmas and definitions that will be  useful throughout this paper.
\medskip

\textbf{Lemma 2.1.} (The reversed H\"{o}lder inequality, \cite{GR}) Let $1<p<\infty$ and $\omega\in A_p(\Bbb{R}^n)$. Then, there exist  positive constants $C$ and $\epsilon$, depending only on $p$ and the $A_p$ condition of $\omega$, such that for any cube $Q$,
$$\left(\frac{1}{|Q|}\int_Q\omega(x)^{1+\epsilon}dx\right)^{\frac{1}{1+\epsilon}}\leq C\left(\frac{1}{|Q|}\int_Q\omega(x)dx\right).\eqno(2.1)$$
\medskip

\textbf{Lemma 2.2.} (\cite{CWX, I})    Let $1\leq p_1, p_2,\cdots, p_m\leq \infty$, $1/p=\sum\limits_{i=1}^m1/p_i$ and $0<q<\infty$. A vector $\vec{\omega}$ of weights satisfies $\vec{\omega}\in A_{(\vec{P},q)}$ if and only if

(i) $\nu_{\vec{\omega}}^q\in A_{1+q(m-\frac{1}{p})}$;

(ii)$\omega_i^{-p_i'}\in A_{1+p_i'\cdot s_i} (i=1,\cdots,m)$
where $s_i=1/q+m-1/p-\frac{1}{p_i'}$.

Moreover, Moen \cite{Mo1} gave another  characterization of $A_{(\vec{P},q)}$.
\medskip

\textbf{Lemma 2.3.} (\cite{Mo1}) Suppose $1<p_1,\cdots,p_m<\infty$ and $\vec{\omega}\in A_{(\vec{P},q)}$. Then
$$\nu_{\vec{\omega}}^q\in A_{mq}\ \ \ \ \text{and} \ \ \ \ \omega^{-p_i'}\in A_{mp_i'}.$$

From Lemmas 2.2 or 2.3, we know that if $[\vec{\omega}]_{aq_0,q,(\frac{sp_1}{s+p_1},\frac{rp_2}{r+p_2})}<\infty$,  then
$$\nu_{\vec{\omega}}^q\in A_{1+q(1-1/p)}\ \ (p>1),\ \ \ \ \  \omega_1^{-r(\frac{p_1}{r})'}\in A_{1+r(\frac{p_1}{r})'(\frac{1}{q}-\frac{1}{p_2}+\frac{1}{s})}, \ \ \ \ \omega_2^{-s(\frac{p_2}{s})'}\in A_{1+s(\frac{p_2}{s})'(\frac{1}{q}-\frac{1}{p_1}+\frac{1}{r})},$$
or
$$\nu_{\vec{\omega}}^q\in A_{2q}\ \ (q>\frac{1}{2}),\ \ \ \ \ \ \ \omega_1^{-r(\frac{p_1}{r})'}\in A_{2r(\frac{p_1}{r})'},\ \ \ \omega_2^{-s(\frac{p_2}{s})'}\in A_{2s(\frac{p_2}{s})'}.$$
Thus, we conclude that the functions $\nu_{\vec{\omega}}^q$, $\omega_1^{-r(\frac{p_1}{r})'}$ and  $\omega_2^{-s(\frac{p_2}{s})'}$ all satisfy the reversed H\"{o}lder inequality throughout the proof of Theorem 1.1.

Next, we introduce some  maximal functions (see \cite{LOPTT} or \cite{MW}).

 The maximal function $M$ and the fractional maximal function $M_\alpha$ are defined by
$$Mf(x)=\sup\limits_{Q\ni x}\frac{1}{|Q|}\int_Q|f(y)|dy$$
and
$$M_\alpha f(x)=\sup\limits_{Q\ni x}\frac{1}{|Q|^{1-\alpha/n}}\int_Q|f(y)|dy, \ \ \ \ (0<\alpha<n)$$
with  $Q$ runs over all cubes containing $x$ respectively.\\
Furthermore, for any $p>1$, we denote
$$M^{(p)}f(x)=\sup\limits_{ Q\ni x}\left(\frac{1}{|Q|}\int_Q|f(y)|^pdy\right)^{1/p}.$$

Before giving the next lemma which is the most important throughout this paper,  we  introduce some notations. First, we define the set of all dyadic grids as follows.

A dyadic grid $\mathcal{D}$ is a countable collection of cubes that satisfies the following properties:

(i) $Q\in \mathcal{D}\Rightarrow l(Q)=2^{-k}$ for some $k\in \Bbb{Z}$.

(ii) For each $k\in \Bbb{Z}$, the set $\{Q\in \mathcal{D}:l(Q)=2^{-k}\}$ forms a partition of $\Bbb{R}^n$.

(iii) $Q,P\in \mathcal{D}\Rightarrow Q\cap P\in \{P,Q,\emptyset\}$.

One very clear example (see \cite{HM,L}) for this concept is the dyadic grid that is formed by translating and then dilating the unit cube $[0,1)^n$ all over $\Bbb{R}^n$. More precisely, it is formulated as
$$\mathcal{D}=\{2^{-k}([0,1)^n+m):k\in \Bbb{Z},m\in \Bbb{Z}^n\}.$$

In practice, we also make extensive use of the following family of dyadic grids.

$$\mathcal{D}^t=\{2^{-k}([0,1)^n+m+(-1)^kt):k\in \Bbb{Z},m\in \Bbb{Z}^n\}, t\in\{0,1/3\}^n.$$

In \cite{L}, Lerner proved the following theorem.
\medskip

\textbf{Lemma 2.4.} (\cite{L}) Given any cube in $\Bbb{R}^n$, there exists a $t\in \{0,1/3\}^n$ and a cube $Q_t\in \mathcal{D}^t$, such that $Q\subset Q_t$ and $l(Q_t)\leq 6l(Q)$.

Next, let us give a decomposition result about cubes.    Suppose that  $Q_0$ is a cube and let $f$ be a function belonging to $L^1_{\text{loc}}(\Bbb{R}^n)$. Then, we set
$$\mathcal{D}(Q_0)\equiv \{Q\in \mathcal{D}: Q\subset Q_0\}.$$

Moreover, we denote that $3Q_0$ is the unique cube concentric to $Q_0$ and having the volume $3^n|Q_0|$. Denote $$m_{3Q_0}(|f|^r,|g|^s)=\left(\frac{1}{|3Q_0|}\int_{3Q_0}|f(x)|^{r}dx\right)^{1/r}\left(\frac{1}{|3Q_0|}\int_{3Q_0}|g(x)|^{s}dx\right)^{1/s},$$
where $p,s>1$ and $1/p+1/s=1$.

Next, we introduce the  sparse family of Calder\'{o}n-Zygmund cubes. More precisely, for each $k\in \Bbb{Z}^+$,
$$D_k\equiv\bigcup \left\{Q: Q\in \mathcal{D}(Q_0), m_{3Q_0}(|f|^r,|g|^s)>a^k\right\},$$
where  $a$ will be chosen later.

Considering the maximal cubes with respect to inclusion, we may write
$$D_k=\bigcup\limits_jQ_{k,j},$$
where the cubes $\{Q_{k,j}\}\subset \mathcal{D}(Q_0)$  are nonoverlapping.  That is, $\{Q_{k,j}\}$ is a family of cubes satisfying
$$\sum\limits_j\chi_{Q_{k,j}}\leq \chi_{Q_0}\eqno(2.2)$$
for almost everywhere. By the maximality of $Q_{k,j}$,  we get
$$
\begin{array}{ll}
a^k<m_{3Q_{k,j}}(|f|^r,|g|^s)<2^{2n}a^k.
\end{array}
\eqno(2.3)
$$
For the properties of $Q_{k,j}$, there is \\
(iv) For any fixed $k$, $Q_{k,j}$ are nonoverlapping for different $j$.\\
(v) If $k_1<k_2$, then there exists $i$, such that $Q_{k_2,j}\subset Q_{k_1,i}$ for any $j\in \Bbb{Z}$.

Next, we will use a clever idea proposed by Tanaka in  \cite{T} to   decompose $Q_0$ as follows.\\

Let $E_0=Q_0\setminus D_1, E_{k,j}=Q_{k,j}\setminus D_{k+1}$. Then,  we have the following lemma.

\medskip

\textbf{Lemma 2.5.} The set $\{E_0\}\bigcup \{E_{k,j}\}$ forms a disjoint family of sets, which decomposes $Q_0$, and satisfies
$$
\begin{array}{ll}
|Q_0|\leq 2|E_0|,  \ \ |Q_{k,j}|\leq 2|E_{k,j}|.
\end{array}
\eqno(2.4)
$$
\textbf{Proof:}  We adopt some basic techniques from \cite{HM} to prove this lemma. By the definitions of $Q_{k,j}$ and $D_{k+1}$, there is
\begin{align*}
&\left|Q_{k,j}\bigcap D_{k+1}\right|=\sum\limits_{Q_{k+1,i}\subset Q_{k,j}}|Q_{k+1,i}|\\
&\leq \frac{1}{a^{k+1}}\sum\limits_{i}\left[\left(|Q_{k+1,i}|\left(\frac{1}{|3Q_{k+1,i}|}\int_{3Q_{k+1,i}}|f(x)|^rdx\right)\right)^{1/r}\right.\\
&\left.\times\left(|Q_{k+1,i}|\left(\frac{1}{|3Q_{k+1,i}|}\int_{3Q_{k+1,i}}|g(x)|^sdx\right)\right)^{1/s}\right]\\
&\leq \frac{1}{a^{k+1}}\left(\sum\limits_{i}|Q_{k+1,i}|\left(\frac{1}{|3Q_{k+1,i}|}\int_{3Q_{k+1,i}}|f(x)|^rdx\right)\right)^{1/r}\\
&\times\left(\sum\limits_{i}|Q_{k+1,i}|\left(\frac{1}{|3Q_{k+1,i}|}\int_{3Q_{k+1,i}}|g(x)|^sdx\right)\right)^{1/s}\\
&\leq \frac{1}{a^{k+1}}\left(|Q_{k,j}|\left(\frac{1}{|3Q_{k,j}|}\int_{3Q_{k,j}}|f(x)|^rdx\right)\right)^{1/r}\left(|Q_{k,j}|\left(\frac{1}{|3Q_{k,j}|}\int_{3Q_{k,j}}|g(x)|^sdx\right)\right)^{1/s},\\
\end{align*}
where the last inequality follows from the fact $Q_{k+1,i}\subset Q_{k,j}$ and $Q_{k,j}$ are nonoverlapping.\\
Then, using (2.3),  we get
\begin{align*}
&\left|Q_{k,j}\bigcap D_{k+1}\right|\\
&\leq \frac{1}{a^{k+1}}\left(|Q_{k,j}|\left(\frac{1}{|3Q_{k,j}|}\int_{3Q_{k,j}}|f(x)|^rdx\right)\right)^{1/r}\left(|Q_{k,j}|\left(\frac{1}{|3Q_{k,j}|}\int_{3Q_{k,j}}|g(x)|^sdx\right)\right)^{1/s}\\
&\leq \frac{2^{2n}}{a^{k+1}}|Q_{k,j}|a^k=\frac{2^{2n}}{a}|Q_{k,j}|.
\end{align*}
Thus, if we choose $a=2^{2n+1}$,   we may get
$$\left|Q_{k,j}\bigcap D_{k+1}\right|\leq \frac{1}{2}|Q_{k,j}|.\eqno(2.5)$$
 Similarly, we can also get
$$\left|D_1\right|\leq \frac{1}{2}|Q_{0}|.\eqno(2.6)$$
Thus, we obtain (2.4) from (2.5) and (2.6).

\medskip

\textbf{Lemma 2.6.} (\cite{A}) Let $0<\alpha<n$, $1<q\leq p<\infty$ and $1<t\leq s<\infty$.  Assume $1/s=1/p-\frac{\alpha}{n}$, $\frac{t}{s}=\frac{q}{p}$. Then, there exists a positive constant $C$ such that
$$\|M_\alpha f\|_{\mathcal{M}^s_t}\leq \|I_\alpha f\|_{\mathcal{M}^s_t}\leq C\|f\|_{\mathcal{M}^p_q}.$$
\medskip

\textbf{Lemma 2.7.} Suppose that there exists real numbers $t, q, p$ satisfying $1<t<q\leq p<\infty$.  Then, we have
$\| f^\ell\|^{1/\ell}_{\mathcal{M}^{p/\ell}_{q/\ell}}=\| f\|_{\mathcal{M}^{p}_{q}}$ with $1<\ell<q$.

\textbf{Proof:} By the definition of Morrey space, we can easily prove Lemma 2.7 and we omit the details here.

\medskip

\textbf{Lemma 2.8.} (\cite{ISST1}) Let $0\leq \alpha<mn$, $\vec{P}=(p_1,\cdots,p_m)$, $\vec{R}=(r_1,\cdots,r_m)$, $0<r_i<p_i<\infty$, $0<q\leq q_0<\infty$, $0<p\leq p_0<\infty$,  $1/q_0=1/p_0-\alpha/n$, $\frac{q}{q_0}=\frac{p}{p_0}$ and $1/p=\sum\limits_{i=1}^m1/p_i$. Then, we have
$$\|\mathcal{M}_{\alpha,\vec{R}}(\vec{f})\|_{\mathcal{M}^{q_0}_{q}}\leq C\|\vec{f}\|_{\mathcal{M}^{p_0}_{\vec{P}}},$$
where
$$\mathcal{M}_{\alpha,\vec{R}}(\vec{f})(x):=\sup\limits_{Q\ni x}l(Q)^\alpha\prod\limits_{i=1}^m\left(\frac{1}{|Q|}\int_Qf_i(y_i)dy_i\right)^{1/r_i}.$$
\section{Proof of Theorem 1.1.}
\noindent For the proof of (1.7), we decompose the proof into two cases: $q>1$ and $q\leq 1$.
\subsection{The case $q>1$.}
Fix a cube $Q_0=Q(x_0,\delta)$ with $\delta>0$. Then, for any $x\in Q_0$, we may decompose $B\mathcal{I}_\alpha$ as
\begin{align*}
B\mathcal{I}_\alpha(f,g)(x)&=\int_{\Bbb{R}^n}\frac{f(x-t)g(x+t)}{|t|^{n-\alpha}}dt\\
&=\int_{|t|\leq 2\delta}\frac{f(x-t)g(x+t)}{|t|^{n-\alpha}}dt+\int_{|t|>2\delta}\frac{f(x-t)g(x+t)}{|t|^{n-\alpha}}dt\\
&=:I+II.
\end{align*}
First, we decompose $II$ as
 \begin{align*}
 II&=\sum\limits_{k=0}^\infty\int_{2\cdot 2^k\delta<|t|\leq 2\cdot2^{k+1}\delta}\frac{f(x-t)g(x+t)}{|t|^{n-\alpha}}dt\\
 &\leq \sum\limits_{k=0}^\infty\frac{1}{(2\cdot2^{k+1}\delta)^{n-\alpha}}\int_{|t|\leq 2\cdot2^{k+1}\delta}|f(x-t)g(x+t)|dt\\
 &\leq \sum\limits_{k=0}^\infty\frac{1}{(2\cdot2^{k+1}\delta)^{n-\alpha}}\left(\int_{|t|\leq 2\cdot2^{k+1}\delta}|f(x-t)|^rdt\right)^{1/r}\left(\int_{|t|\leq 2\cdot2^{k+1}\delta}|g(x+t)|^sdt\right)^{1/s}.\\
 \end{align*}
Then, by a change of variables and the fact $x\in Q_0=Q(x_0,\delta)$, we obtain
\begin{align*}
&|Q_0|^{\frac{1}{q_0}-\frac{1}{q}}\left(\int_{Q_0}\left|II\right|^q\left(\prod\limits_{i=1}^2\omega_i(x)\right)^qdx\right)^{1/q}\\
&\leq C|Q_0|^{\frac{1}{q_0}-\frac{1}{q}}\sum\limits_{k=0}^\infty(2\cdot2^{k}\delta)^{\alpha-n}\left(\int_{Q_0}\prod\limits_{i=1}^2\omega_i(x)^qdx\right)^{1/q}\\
&\times \left(\int_{2^{k+3} Q_0}|f(u)|^rdu\right)^{1/r}\left(\int_{2^{k+3} Q_0}|g(v)|^sdv\right)^{1/s}.
\end{align*}
 For $\left(\int_{2^{k+3} Q_0}|f(u)|^rdu\right)^{1/r}$, by the H\"{o}lder inequality,  there is
 \begin{align*}
 &\left(\int_{2^{k+3} Q_0}|f(u)|^rdu\right)^{1/r}\\
 &\leq \left(\int_{2^{k+3} Q_0}|f(u)\omega(u)|^{p_1}du\right)^{1/p_1}\left(\int_{2^{k+3} Q_0}|\omega(u)|^{-r(\frac{p_1}{r})'}du\right)^{1/r-1/p_1}.
 \end{align*}
 Similarly, we have
 \begin{align*}
 &\left(\int_{2^{k+3} Q_0}|g(v)|^sdu\right)^{1/s}\\
 &\leq \left(\int_{2^{k+3} Q_0}|g(v)\omega(v)|^{p_2}dv\right)^{1/p_2}\left(\int_{2^{k+3} Q_0}|\omega(v)|^{-s(\frac{p_2}{s})'}dv\right)^{1/s-1/p_2}.
 \end{align*}
Thus, using the condition $[\vec{\omega}]_{aq_0,q,(\frac{sp_1}{s+p_1},\frac{rp_2}{r+p_2})}<\infty$ and $a>1$, we get
\begin{align*}
&|Q_0|^{\frac{1}{q_0}-\frac{1}{q}}\left(\int_{Q_0}\left|II\right|^q\left(\prod\limits_{i=1}^2\omega_i(x)\right)^qdx\right)^{1/q}\\
&\leq C\sum\limits_{k=1}^\infty(2^{k+2}\delta)^{\alpha-n}|Q_0|^{1/q_0-1/q}\left(\int_{Q_0}\prod\limits_{i=1}^2\omega_i(x)^qdx\right)^{1/q}\\
&\times \left(\int_{2^{k+3} Q_0}|f(u)\omega(u)|^{p_1}du\right)^{1/p_1}\left(\int_{2^{k+3} Q_0}|\omega(u)|^{-r(\frac{p_1}{r})'}du\right)^{1/r-1/p_1}\\
&\times \left(\int_{2^{k+3} Q_0}|g(v)\omega(v)|^{p_2}dv\right)^{1/p_2}\left(\int_{2^{k+3} Q_0}|\omega(v)|^{-s(\frac{p_2}{s})'}dv\right)^{1/s-1/p_2}\\
&\leq C\|(f_1\omega_1,f_2\omega_2)\|_{\mathcal{M}^{p_0}_{\vec{P}}}\sum\limits_{k=1}^\infty(2^k\delta)^{\alpha-n}|Q_0|^{1/q_0-1/q+1/q}|2^{k+3} Q_0|^{1/p-1/p_0+1/r-1/p_1+1/s-1/p_2}\\
&\times \left(\frac{|Q_0|}{|2^{k+3} Q_0|}\right)^{-\frac{1}{aq_0}}\left(\frac{|Q_0|}{|2^{k+3} Q_0|}\right)^{\frac{1}{aq_0}}\left(\frac{1}{|Q_0|}\int_{Q_0}\prod\limits_{i=1}^2\omega_i(x)^q\right)^{1/q}\\
&\times \left(\frac{1}{|2^{k+3} Q|}\int_{2^{k+3} Q_0}|f(u)\omega(u)|^{p_1}du\right)^{1/p_1}\left(\int_{2^{k+3} Q_0}|\omega(u)|^{-r(\frac{p_1}{r})'}du\right)^{1/r-1/p_1}\\
&\times \left(\frac{1}{|2^{k+3} Q|}\int_{2^{k+3} Q_0}|g(v)\omega(v)|^{p_2}dv\right)^{1/p_2}\left(\int_{2^{k+3} Q_0}|\omega(v)|^{-s(\frac{p_2}{s})'}dv\right)^{1/s-1/p_2}\\
&\leq C[\vec{\omega}]_{aq_0,q,(\frac{sp_1}{s+p_1},\frac{rp_2}{r+p_2})}\|(f_1\omega_1,f_2\omega_2)\|_{\mathcal{M}^{p_0}_{\vec{P}}},
\end{align*}
which implies
$$\|II\cdot\nu_{\vec{\omega}}\|_{\mathcal{M}^{q_0}_q}\leq C[\vec{\omega}]_{aq_0,q,(\frac{sp_1}{s+p_1},\frac{rp_2}{r+p_2})}\|(f_1\omega_1,f_2\omega_2)\|_{\mathcal{M}^{p_0}_{\vec{P}}}.\eqno(3.1)$$

It remains to give the  estimates of $\|I\cdot\nu_{\vec{\omega}}\|_{\mathcal{M}^{q_0}_q}$. First, we prove  the following lemma.
\medskip

\textbf{Lemma 3.1.} Denote  $I=\int_{|t|\leq 2\delta}\frac{f(x-t)g(x+t)}{|t|^{n-\alpha}}dt$  and $Q_0=Q(x_0,\delta)$ with $\delta>0$. There exists a positive constant independent of $f$ and $g$, such that
 $$I\leq C\sum\limits_{Q\in \mathcal{D}(Q_0)}l(Q)^\alpha m_{3Q}(|f|^r,|g|^s)\chi_Q(x).\eqno(3.2)$$

\textbf{Proof.} By the definition $I$, we may get

\begin{align*}
I&=\int_{|t|\leq 2\delta}\frac{f(x-t)g(x+t)}{|t|^{n-\alpha}}dt\\
&=\sum\limits_{k=0}^{+\infty}\int_{2\cdot 2^{-k-1}\delta<|t|\leq 2\cdot 2^{-k}\delta}\frac{f(x-t)g(x+t)}{|t|^{n-\alpha}}dt\\
&\leq \sum\limits_{k=0}^{+\infty}\frac{1}{(2\cdot2^k\delta)^{n-\alpha}}\int_{|t|\leq 2\cdot2^{-k}\delta}f(x-t)g(x+t)dt\\
&\leq  \sum\limits_{k=0}^{+\infty}\frac{1}{(2\cdot2^{-k}\delta)^{n-\alpha}}\left(\int_{|t|\leq 2\cdot 2^{-k}\delta}|f(x-t)|^rdt\right)^{1/r}\left(\int_{|t|\leq 2\cdot 2^{-k}\delta}|g(x+t)|^sdt\right)^{1/s}.
\end{align*}
Then, by a change of variables and the fact $x\in Q_0$, it is easy to see
\begin{align*}
I&\leq C\sum\limits_{k=0}^{+\infty}\frac{1}{(2\cdot2^{-k}\delta)^{n-\alpha}}\left(\int_{|u-x|\leq 2\cdot2^{-k}\delta}|f(u)|^rdu\right)^{1/r}\left(\int_{|v-x|\leq 2\cdot2^{-k}\delta}|g(v)|^sdv\right)^{1/s}\\
&\leq C\sum\limits_{k=0}^{+\infty}\sum\limits_{\substack{Q\in \mathcal{D}(Q_0)\\l(Q)=2^{-k}\delta}}l(Q)^{\alpha-n}\left(\int_{|u-x|\leq 2l(Q)}|f(u)|^rdu\right)^{1/r}\\
&\times\left(\int_{|v-x|\leq 2l(Q)}|g(v)|^sdv\right)^{1/s}\chi_Q(x)\\
&\leq C\sum\limits_{k=0}^{+\infty}\sum\limits_{\substack{Q\in \mathcal{D}(Q_0)\\l(Q)=2^{-k}\delta}}l(Q)^{\alpha-n}\left(\int_{3Q}|f(u)|^rdu\right)^{1/r}\left(\int_{3Q}|g(v)|^sdv\right)^{1/s}\chi_Q(x)\\
&= C\sum\limits_{Q\in \mathcal{D}(Q_0)}l(Q)^\alpha m_{3Q}(|f|^r,|g|^s)\chi_Q(x),
\end{align*}
which implies
$$I\leq C\sum\limits_{Q\in \mathcal{D}(Q_0)}l(Q)^\alpha m_{3Q}(|f|^r,|g|^s)\chi_Q(x).$$
Thus,  the proof of Lemma 3.1 has been finished.
\medskip

Next,  we recall some notations from Section 2.  For $r,s>1$ with $1/r+1/s=1$, we set
$$\mathcal{D}_0(Q_0)\equiv \{Q\in \mathcal{D}(Q_0): m_{3Q}(|f|^r,|g|^s)\leq a\}$$
and
$$\mathcal{D}_{k,j}(Q_0)\equiv \{Q\in \mathcal{D}(Q_0): Q\subset Q_{k,j}, a^{k}<m_{3Q}(|f|^r,|g|^s)\leq a^{k+1}\},$$
where $a$ is the same as in Section 2.
Thus, we have
$$\mathcal{D}(Q_0)=\mathcal{D}_0(Q_0)\bigcup \left(\bigcup\limits_{k,j}\mathcal{D}_{k,j}(Q_0)\right).$$
As $q>1$, by duality, there is
$$\left(\int_{Q_0}|I|^q(\omega_1(x)\omega_2(x))^qdx\right)^{1/q}=\sup\limits_{\|h\|_{L^{q'}(Q_0)\leq 1}}\|I\omega_1\omega_2h\|_{L^1(Q_0)}.$$
Then, we denote
$$I_0:=\sum\limits_{Q\in \mathcal{D}_0(Q_0)}l(Q)^\alpha m_{3Q}(|f|^r,|g|^s)\int_Q\omega_1(x)\omega_2(x)h(x)dx,$$
and
$$I_{k,j}:=\sum\limits_{Q\in \mathcal{D}_{k,j}(Q_0)}l(Q)^\alpha m_{3Q}(|f|^r,|g|^s)\int_Q\omega_1(x)\omega_2(x)h(x)dx.$$
From $\mathcal{D}(Q_0)=\mathcal{D}_0(Q_0)\bigcup \left(\bigcup\limits_{k,j}\mathcal{D}_{k,j}(Q_0)\right)$ and (3.2), we get
$$\left(\int_{Q_0}|I|^q|\omega_1(x)\omega_2(x)|^qdx\right)^{1/q}\leq I_0+I_{k,j}.\eqno(3.3)$$
For $I_{k,j}$, recall that $q>1$, $a>1$ and $\alpha>0$. Then, using (2.3),  the H\"{o}lder inequality, Lemmas 2.5 and the property of $\mathcal{D}$,  we obtain
\begin{align*}
I_{k,j}&\leq a^{k+1}\sum\limits_{Q\in \mathcal{D}_{k,j}(Q_0)}l(Q)^{\alpha}\int_Q\omega_1(x)\omega_2(x)h(x)dx\\
&\leq Ca^{k+1}l(Q_{k,j})^{\alpha}\int_{Q_{k,j}}\omega_1(x)\omega_2(x)h(x)dx\\
&\leq Cam_{3Q_{k,j}}(|f|^r,|g|^s)l(Q_{k,j})^\alpha\int_{Q_{k,j}}\omega_1(x)\omega_2(x)h(x)dx\\
&\leq Cam_{3Q_{k,j}}(|f|^r,|g|^s)l(Q_{k,j})^\alpha|Q_{k,j}|\\
&\times\left(\frac{1}{|Q_{k,j}|}\int_{Q_{k,j}}(\omega_1(x)\omega_2(x))^{aq}dx\right)^{1/aq}\left(\frac{1}{|Q_{k,j}|}\int_{Q_{k,j}}|h(x)|^{(aq)'}dx\right)^{1/(aq)'}\\
&\leq Ca|E_{k,j}|m_{3Q_{k,j}}(|f|^r,|g|^s)l(Q_{k,j})^\alpha\\
&\times\left(\frac{1}{|Q_{k,j}|}\int_{Q_{k,j}}(\omega_1(x)\omega_2(x))^{aq}dx\right)^{1/aq}\left(\frac{1}{|Q_{k,j}|}\int_{Q_{k,j}}|h(x)|^{(aq)'}dx\right)^{1/(aq)'}\\
\end{align*}
\begin{align*}
&\leq Ca\int_{E_{k,j}}m_{3Q_{k,j}}(|f|^r,|g|^s)l(Q_{k,j})^\alpha\left(\frac{1}{|Q_{k,j}|}\int_{Q_{k,j}}(\omega_1(y)\omega_2(y))^{aq}dy\right)^{1/aq}\\
&\times \left(\frac{1}{|Q_{k,j}|}\int_{Q_{k,j}}|h(y)|^{(aq)'}dy\right)^{1/(aq)'}dx\\
&\leq Ca\int_{E_{k,j}}\left[M(h^{(aq)'})(x)\right]^{\frac{1}{(aq)'}}\tilde{M}^{aq}_{\alpha,r,s}(\vec{f},\vec{\omega})(x)dx,
\end{align*}
where
$$\tilde{M}^{aq}_{\alpha,r,s}(\vec{f},\vec{\omega})(x)=\sup\limits_{Q\ni x}l(Q)^\alpha m_{3Q}(|f|^r,|g|^s)\left(\frac{1}{|Q|}\int_Q(\omega_1(y)\omega_2(y))^{aq}dy\right)^{1/aq}.$$
Similarly, there is
$$I_0\leq Ca\int_{E_{0}}\left[M(h^{(aq)'})(x)\right]^{\frac{1}{(aq)'}}\tilde{M}^{aq}_{\alpha,r,s}(\vec{f},\vec{\omega})(x)dx.$$
Thus, by the boundedness of the Hardy-Littlewood maximal function, the H\"{o}lder inequality and the fact $q'>(aq)'$, we obtain
$$
\begin{array}{ll}
&I_0+\sum\limits_{k,j}I_{k,j}\\
&\leq C\left(\int_{Q_0}|\tilde{M}^{aq}_{\alpha,r,s}(\vec{f},\vec{\omega})(x)|^qdx\right)^{1/q}\left(\int_{Q_0}\left[M(h^{(aq)'})(x)\right]^{\frac{q'}{(aq)'}}\right)^{1/q'}\\
&\leq C\left(\int_{Q_0}|\tilde{M}^{aq}_{\alpha,r,s}(\vec{f},\vec{\omega})(x)|^qdx\right)^{1/q}\left(\int_{Q_0}|h(x)|^{q'}dx\right)^{1/q'}\\
&\leq C\|\tilde{M}^{aq}_{\alpha,r,s}(\vec{f},\vec{\omega})\|_{L^{q}(Q_0)}.
\end{array}
\eqno(3.4)
$$
Using the H\"{o}lder inequality and the reversed H\"{o}lder inequality for  $\omega_1^{-r(\frac{p_1}{r})'}$ and  $\omega_2^{-s(\frac{p_2}{s})'}$, we have
\begin{align*}
&m_{3Q}(|f|^r,|g|^s)=\left(\frac{1}{|3Q|}\int_{3Q}|f|^rdx\right)^{1/r}\left(\frac{1}{|3Q|}\int_{3Q}|g|^sdx\right)^{1/s}\\
&\leq |3Q|^{-1}\left(\int_{3Q}\left(|f(x)|^r\omega_1(x)^r\right)^{\frac{p_1r}{ar}}dx\right)^{\frac{1}{r}\frac{ar}{p_1}}\left(\int_{3Q}\omega_1(x)^{-r\frac{p_1}{p_1-ar}}dx\right)^{\frac{1}{r}-\frac{a}{p_1}}\\
&\times\left(\int_{3Q}\left(|g(x)|^s\omega_2(x)^s\right)^{\frac{p_2s}{as}}dx\right)^{\frac{1}{s}\frac{as}{p_2}}\left(\int_{3Q}\omega_2(x)^{-s\frac{p_2}{p_2-as}}dx\right)^{\frac{1}{s}-\frac{a}{p_2}}\\
&\leq |3Q|^{-1}|3Q|^{a/p_1+a/p_2+1/r-a/p_1+1/s-a/p_2}\left(\frac{1}{|3Q|}\int_{3Q}\left(|f(x)|^r\omega_1(x)^r\right)^{\frac{p_1}{ar}r}dx\right)^{\frac{1}{r}\frac{ar}{p_1}}\\
&\times\left(\frac{1}{|3Q|}\int_{3Q}\left(|g(x)|^s\omega_2(x)^s\right)^{\frac{p_2}{as}s}dx\right)^{\frac{1}{s}\frac{as}{p_2}}\left(\frac{1}{|3Q|}\int_{3Q}\omega_2(x)^{-s\frac{p_2}{p_2-as}}dx\right)^{\frac{1}{s}-\frac{a}{p_2}}\\
&\times\left(\frac{1}{|3Q|}\int_{3Q}\omega_1(x)^{-r\frac{p_1}{p_1-ar}}dx\right)^{\frac{1}{r}-\frac{a}{p_1}}\\
\end{align*}
\begin{align*}
&\leq \left(\frac{1}{|3Q|}\int_{3Q}|f(x)\omega_1(x)|^{\frac{p_1}{a}}dx\right)^{\frac{a}{p_1}}\left(\frac{1}{|3Q|}\int_{3Q}|g(x)\omega_2(x)|^{\frac{p_2}{a}}dx\right)^{\frac{a}{p_2}}\\
&\times \left(\frac{1}{|3Q|}\int_{3Q}\omega_1(x)^{-\frac{rp_1}{p_1-r}}dx\right)^{\frac{1}{r}-\frac{1}{p_1}}\left(\frac{1}{|3Q|}\int_{3Q}\omega_2(x)^{-\frac{sp_2}{p_2-s}}dx\right)^{\frac{1}{s}-\frac{1}{p_2}}.
\end{align*}
Recall the definition of $\mathcal{M}_{\alpha,\vec{R}}(\vec{f})(x)$ in Section 2, we obtain
$$\mathcal{M}^{aq}_{\alpha,r,s}(\vec{f},\vec{\omega})(x)\leq C[\vec{\omega}]_{aq_0,q,(\frac{sp_1}{s+p_1},\frac{rp_2}{r+p_2})}\mathcal{M}_{\alpha,\frac{\vec{P}}{a}}(f_1\omega_1,f_2\omega_2)(x).\eqno(3.5)$$
Using Lemma 2.8 and (3.3)-(3.5), we have
$$
\begin{array}{ll}
\|I\cdot\nu_{\vec{\omega}}\|_{\mathcal{M}^{q_0}_q}&\leq C[\vec{\omega}]_{aq_0,q,(\frac{sp_1}{s+p_1},\frac{rp_2}{r+p_2})}\|\mathcal{M}_{\alpha,\frac{\vec{P}}{a}}(f_1\omega_1,f_2\omega_2)\|_{\mathcal{M}^{q_0}_q}\\
&\leq C[\vec{\omega}]_{aq_0,q,(\frac{sp_1}{s+p_1},\frac{rp_2}{r+p_2})}\|(f_1\omega_1,f_2\omega_2)\|_{\mathcal{M}_{\vec{P}}^{p_0}}.
\end{array}
\eqno(3.6)
$$
Combining (3.1) and (3.6), we finish the proof of Theorem 1.1 for the case $q>1$.
\subsection{The case $q\leq 1$.}
First, we denote
$$L:=\left(\sum\limits_{Q\in \mathcal{D}(Q_0)}l(Q)^\alpha m_{3Q}(|f|^r,|g|^s)\chi_Q(x)\right)^q.$$
Since $q\leq 1$, there is
\begin{align*}
L&\leq \sum\limits_{Q\in \mathcal{D}(Q_0)}l(Q)^{q\alpha}m_{3Q}(|f|^r,|g|^s)^q\chi_Q(x)\\
&\leq \left(\sum\limits_{Q\in \mathcal{D}_0(Q_0)}+\sum\limits_{k,j}\sum\limits_{Q\in\mathcal{D}_{k,j}(Q_0)}\right)l(Q)^{q\alpha}m_{3Q}(|f|^r,|g|^s)^q\chi_Q(x).
\end{align*}
Recall that $\nu_{\vec{\omega}}(x)=\omega_1(x)\omega_2(x)$. Then,  we obtain
\begin{align*}
&\int_{Q_0}|B\mathcal{I}_\alpha(f,g)(x)|^q(\omega_1(x)\omega_2(x))^qdx\\
&\leq C\left(\sum\limits_{Q\in \mathcal{D}_0(Q_0)}+\sum\limits_{k,j}\sum\limits_{Q\in\mathcal{D}_{k,j}(Q_0)}\right)l(Q)^{\alpha q}m_{3Q}(|f|^r,|g|^s)^q\int_{Q}(\omega_1(x)\omega_2(x))^qdx\\
&:=C(I_0'+\sum\limits_{k,j}I_{k,j}').
\end{align*}
For $I_{k,j}'$, there is
\begin{align*}
I_{k,j}'&=\sum\limits_{Q\in\mathcal{D}_{k,j}(Q_0)}l(Q)^{\alpha q}m_{3Q}(|f|^r, |g|^s)^q\int_{Q}(\omega_1(x)\omega_2(x))^qdx\\
&\leq l(Q_{k,j})^{\alpha q}(a^{k+1})^q\int_{Q_{k,j}}(\omega_1(x)\omega_2(x))^qdx\\
&\leq Ca|Q_{k,j}|l(Q_{k,j})^{\alpha q}(a^{k+1})^qm_{3Q_{k,j}}(|f|^r,|g|^s)^q\left(\frac{1}{|Q_{k,j}|}\int_{Q_{k,j}}(\omega_1(x)\omega_2(x))^qdx\right)\\
\end{align*}
\begin{align*}
&\leq Ca|E_{k,j}|l(Q_{k,j})^{\alpha q}m_{3Q_{k,j}}(|f|^r,|g|^s)^q\left(\frac{1}{|Q_{k,j}|}\int_{Q_{k,j}}(\omega_1(x)\omega_2(x))^qdx\right)\\
&\leq Ca\int_{E_{k,j}}\left[l(Q_{k,j})^\alpha m_{3Q_{k,j}}(|f|^r,|g|^s)\left(\frac{1}{|Q_{k,j}|}\int_{Q_{k,j}}(\omega_1(y)\omega_2(y))^qdy\right)^{1/q}\right]^qdx\\
&\leq Ca\int_{E_{k,j}}\tilde{\mathcal{M}}^q_{\alpha,r,s}(\vec{f},\omega_1,\omega_2)(x)^qdx,
\end{align*}
where
$$\tilde{\mathcal{M}}^q_{\alpha,r,s}(\vec{f},\omega_1,\omega_2)(x)=\sup\limits_{Q\ni x}l(Q)^\alpha m_{3Q}(|f|^r,|g|^s)\left(\frac{1}{|Q|}\int_Q(\omega_1(y)\omega_2(y))^qdy\right)^{1/q}.$$
Similarly, there is
$$I_0'\leq Ca\int_{E_0}\tilde{\mathcal{M}}^q_{\alpha,r,s}(\vec{f},\vec{\omega})(x)^qdx.$$
Thus,  we obtain
$$I_0'+\sum\limits_{k,j}I_{k,j}'\leq C\int_{Q_0}\tilde{\mathcal{M}}^q_{\alpha,r,s}(\vec{f},\vec{\omega})(x)^qdx.$$
Then, by a similar argument as in the proof of (3.5), there is
$$\tilde{\mathcal{M}}^q_{\alpha,r,s}(\vec{f},\vec{\omega})(x)\leq [\vec{\omega}]_{aq_0,q,(\frac{sp_1}{s+p_1},\frac{rp_2}{r+p_2})}\mathcal{M}_{\alpha,\frac{\vec{P}}{a}}(f_1\omega_1,f_2\omega_2)(x).\eqno(3.7)$$
Now, using Lemma 2.8 and the definition of the Morrey space, we finish the proof of Theorem 1.1 with $q\leq 1$.
\section{Two-weight norm inequalities for $B\mathcal{I}_\alpha$.}
In this section, we are going to give the two-weight norm inequalities for $B\mathcal{I}_\alpha$ on the  Morrey type spaces. Suppose that $v$ and $\vec{\omega}=(\omega_1,\omega_2)$ satisfy the following condition:
$$[v,\vec{\omega}]_{q_0,q,\vec{P}}:=\sup\limits_{\substack{Q\subset Q'\\Q,Q':\text{cubes}}}\left(\frac{|Q|}{|Q'|}\right)^{1/q_0}\left(\frac{1}{|Q|}\int_Qv(x)^qdx\right)^{1/q}\prod\limits_{i=1}^2\left(\frac{1}{|Q|}\int_{Q'}\omega_i(y_i)^{-p_i'}dy_i\right)^{1/p_i'}.$$
Obviously, if  $[v,\vec{\omega}]_{q_0,q,(\frac{sp_1}{s+p_1},\frac{rp_2}{r+p_2})}<\infty$, we cannot get the reversed H\"{o}lder inequality for $v$, $\omega_1^{-r(\frac{p_1}{r})'}$ and  $\omega_2^{-s(\frac{p_2}{s})'}$.\\ Thus, by checking the proof of Theorem 1.1, we obtain

\textbf{Theorem 4.1.} Suppose $0<\alpha<n$, $p_1>r>1$, $p_2>s>1$, $1/r+1/s=1$, $1/p=1/p_1+1/p_2$, $1<p_1,p_2<\infty$,  $0<p\leq p_0<\infty$, $0<q\leq q_0<\infty$. Assume that
$$1/q_0=1/p_0-\alpha/n, \ \ q/q_0=p/p_0.$$

\textbf{Case 1.} If $q>1$, suppose that  there exists $a$ satisfying $1<a<\min\{\frac{p_1}{s'},\frac{p_2}{r'}\}$, such that
$$[v,\vec{\omega}]_{aq_0,aq,(\frac{sp_1}{as+p_1},\frac{rp_2}{ar+p_2})}<\infty.$$
Then, there exists a positive constant $C$ independent of $f$ and $g$, such that
$$\|B\mathcal{I}_{\alpha}(f,g)v\|_{\mathcal{M}^{q_0}_{q}}\leq C[v,\vec{\omega}]_{aq_0,aq,(\frac{sp_1}{as+p_1},\frac{rp_2}{ar+p_2})}\|(f\omega_1,g\omega_2)\|_{\mathcal{M}^{p_0}_{\vec{P}}}.\eqno(4.1)$$

\textbf{Case 2.}  If $0<q\leq 1$, suppose that  there exists $a$ satisfying $1<a<\min\{\frac{p_1}{s'},\frac{p_2}{r'}\}$, such that
$$[v,\vec{\omega}]_{aq_0,q,(\frac{sp_1}{as+p_1},\frac{rp_2}{ar+p_2})}<\infty.$$
Then, there exists a positive constant $C$ independent of $f$ and $g$, such that
$$\|B\mathcal{I}_{\alpha}(f,g)v\|_{\mathcal{M}^{q_0}_{q}}\leq C[v,\vec{\omega}]_{aq_0,q,(\frac{sp_1}{as+p_1},\frac{rp_2}{ar+p_2})}\|(f\omega_1,g\omega_2)\|_{\mathcal{M}^{p_0}_{\vec{P}}}.\eqno(4.2)$$

In order to prove Theorem 4.1, recall the definition of $\tilde{M}^{aq}_{\alpha,r,s}(\vec{f},\vec{\omega})(x)$ in Section 3, we need the following lemma.
\medskip

\textbf{Lemma 4.1.} Under the same conditions as in  Theorem 4.1, we have the following estimates for $\tilde{M}^{aq}_{\alpha,r,s}$.

\textbf{Case 1.} For the case $q>1$, suppose that  there exists $a$ satisfying $1<a<\min\{\frac{p_1}{s'},\frac{p_2}{r'}\}$, such that
$$[v,\vec{\omega}]_{aq_0,aq,(\frac{sp_1}{as+p_1},\frac{rp_2}{ar+p_2})}<\infty.$$
Then
$$\tilde{M}^{aq}_{\alpha,r,s}(\vec{f},\vec{\omega})(x)\leq C[v,\vec{\omega}]_{aq_0,aq,(\frac{sp_1}{as+p_1},\frac{rp_2}{ar+p_2})}\mathcal{M}_{\alpha,\frac{\vec{P}}{a}}(f_1\omega_1,f_2\omega_2)(x).$$

\textbf{Case 2.} For the case $q\leq 1$, suppose that  there exists $a$ satisfying $1<a<\min\{\frac{p_1}{s'},\frac{p_2}{r'}\}$, such that
$$[v,\vec{\omega}]_{aq_0,q,(\frac{sp_1}{as+p_1},\frac{rp_2}{ar+p_2})}<\infty.$$
Then
$$\tilde{M}^{q}_{\alpha,r,s}(\vec{f},\vec{\omega})(x)\leq C[v,\vec{\omega}]_{aq_0,q,(\frac{sp_1}{as+p_1},\frac{rp_2}{ar+p_2})}\mathcal{M}_{\alpha,\frac{\vec{P}}{a}}(f_1\omega_1,f_2\omega_2)(x).$$
If we check the proof of (3.5) and (3.7) carefully, we can easily get Lemma 4.1 and we omit the details here.\\
Moreover, we can generalize Theorem 4.1 to a more general case.

Suppose that another quantity of two-weight type multiple weights $[v,\vec{\omega}]_{q_0,r_0,q,\vec{P}}$:
$$[v,\vec{\omega}]_{q_0,r_0,q,\vec{P}}:=\sup\limits_{\substack{Q\subset Q'\\Q,Q':\text{cubes}}}\left(\frac{|Q|}{|Q'|}\right)^{1/q_0}|Q'|^{1/r_0}\left(\frac{1}{|Q|}\int_Qv(x)^qdx\right)^{1/q}\prod\limits_{i=1}^2\left(\frac{1}{|Q|}\int_{Q'}\omega_i(y_i)^{-p_i'}dy_i\right)^{1/p_i'}.$$
By checking the proof of Theorem 1.1 again, we have

\medskip

\textbf{Theorem 4.2.} Suppose $0<\alpha<n$, $p_1>r>1$, $p_2>s>1$, $1/r+1/s=1$, $1/p=1/p_1+1/p_2$, $1<p_1,p_2<\infty$,  $0<p\leq p_0<\infty$, $0<q\leq q_0<\infty$. Assume that
$$q/q_0=p/p_0,\ \ 1/q_0=1/p_0+1/r_0-\alpha/n,   \ \ r_0\geq \frac{n}{\alpha}.$$

\textbf{Case 1.} If $q>1$, suppose that  there exists $a$ satisfying $1<a<\min\{\frac{r_0}{q_0}, \frac{p_1}{s'},\frac{p_2}{r'}\}$, such that
$$[v,\vec{\omega}]_{aq_0,r_0,aq,(\frac{sp_1}{as+p_1},\frac{rp_2}{ar+p_2})}<\infty.$$
Then, there exists a positive constant $C$ independent of $f$ and $g$, such that
$$\|B\mathcal{I}_{\alpha}(f,g)v\|_{\mathcal{M}^{q_0}_{q}}\leq C[v,\vec{\omega}]_{aq_0,r_0,aq,(\frac{sp_1}{as+p_1},\frac{rp_2}{ar+p_2})}\|(f\omega_1,g\omega_2)\|_{\mathcal{M}^{p_0}_{\vec{P}}}.\eqno(4.3)$$

\textbf{Case 2.}  If $0<q\leq 1$, suppose that  there exists $a$ satisfying $1<a<\min\{\frac{r_0}{q_0},\frac{p_1}{s'},\frac{p_2}{r'}\}$, such that
$$[v,\vec{\omega}]_{aq_0,r_0,q,(\frac{sp_1}{as+p_1},\frac{rp_2}{ar+p_2})}<\infty.$$
Then, there exists a positive constant $C$ independent of $f$ and $g$, such that
$$\|B\mathcal{I}_{\alpha}(f,g)v\|_{\mathcal{M}^{q_0}_{q}}\leq C[v,\vec{\omega}]_{aq_0,r_0, q,(\frac{sp_1}{as+p_1},\frac{rp_2}{ar+p_2})}\|(f\omega_1,g\omega_2)\|_{\mathcal{M}^{p_0}_{\vec{P}}}.\eqno(4.4)$$
Similarly, to prove Theorem 4.2, we need the following lemma.
\medskip

\textbf{Lemma 4.2.} Under the same conditions as in Theorem 4.2, we have the following estimates for $\tilde{M}^{aq}_{\alpha,r,s}$.

\textbf{Case 1.} For the case $q>1$, suppose that  there exists $a$ satisfying $1<a<\min\{\frac{r_0}{q_0}, \frac{p_1}{s'},\frac{p_2}{r'}\}$, such that
$$[v,\vec{\omega}]_{aq_0,r_0,aq,(\frac{sp_1}{as+p_1},\frac{rp_2}{ar+p_2})}<\infty.$$
Then
$$\tilde{M}^{aq}_{\alpha,r,s}(\vec{f},\vec{\omega})(x)\leq C[v,\vec{\omega}]_{aq_0,r_0,aq,(\frac{sp_1}{as+p_1},\frac{rp_2}{ar+p_2})}\mathcal{M}_{\alpha-\frac{n}{r_0},\frac{\vec{P}}{a}}(f_1\omega_1,f_2\omega_2)(x).$$

\textbf{Case 2.} For the case $0<q\leq 1$, suppose that  there exists $a$ satisfying $1<a<\min\{\frac{r_0}{q_0}, \frac{p_1}{s'},\frac{p_2}{r'}\}$, such that
$$[v,\vec{\omega}]_{aq_0,r_0,q,(\frac{sp_1}{as+p_1},\frac{rp_2}{ar+p_2})}<\infty.$$
Then
$$\tilde{M}^{q}_{\alpha,r,s}(\vec{f},\vec{\omega})(x)\leq C[v,\vec{\omega}]_{aq_0,r_0,q,(\frac{sp_1}{as+p_1},\frac{rp_2}{ar+p_2})}\mathcal{M}_{\alpha-\frac{n}{r_0},\frac{\vec{P}}{a}}(f_1\omega_1,f_2\omega_2)(x).$$

\textbf{Remark 4.3.} For the case $0<q\leq 1$, the results of (4.2) and (4.4) are still different from \cite[Theorem 4.2]{HY}.

\section{ An Olsen type inequality for $B\mathcal{I}_\alpha$.}
In this section, we will give an Olsen type inequality for $B\mathcal{I}_\alpha$. Recall the fractional integral
$$
\begin{array}{ll}
I_{\alpha}(f)(x)=\int_{\Bbb{R}^n}\frac{f(y)}{|x-y|^{n-\alpha}}dy,  \ \ \ \ \ 0<\alpha<n.
\end{array}
\eqno(5.1)
$$

For the study of $I_\alpha$ on the Morrey space, one may see \cite{A, IKS1, IKS2} et al. to find more details.

Particularly, recently  Sawano, Sugano and Tanaka obtained the following
result.

\medskip

\textbf{Theorem A.}     (\cite{SST1}) Suppose that the indices \ $\alpha
,p_{0},q_{0},r_{0},p,q,r_1$ \ satisfy%
\[
1<p\leq p_{0}<\infty ,\text{ }1<q\leq q_{0}<\infty ,1<r_1\leq r_{0}<\infty
\]%
and
\[
r_1>q,\text{ }1/p_{0}>\alpha /n\geq 1/r_{0}.
\]%
Also assume
\[
q/q_{0}=p/p_{0},\text{ }1/p_{0}+1/r_{0}-\alpha /n=1/q_{0}.
\]%
Then, for all \ $f\in \mathcal{M}_{p}^{p_{0}}(%
\mathbb{R}
^{n})$ and \ $h\in \mathcal{M}_{r_1}^{r_{0}}(%
\mathbb{R}
^{n}),$ \
\[
\left\Vert h\cdot I_{\alpha }(f)\right\Vert _{\mathcal{M}_{q}^{q_{0}}(%
\mathbb{R}
^{n})}\leq C\left\Vert f\right\Vert _{\mathcal{M}_{p}^{p_{0}}(%
\mathbb{R}
^{n})}\left\Vert h\right\Vert _{\mathcal{M}_{r_1}^{r_{0}}(%
\mathbb{R}
^{n})},\eqno(5.2)
\]%
where $C$ is a positive constant independet of $f$ and $g$.

The above inequality is called an inequality of Olsen type, since it is
initially proposed by Olsen in \cite{O} and Olsen found that this inequality
would perform crucially in the study of Schr\"{o}dinger equation.  The
inequality was proved in the case \ $n=3$ \ by Conlon and Redondo in \cite{CR} essentially. In fact,
an analogous inequality on a generalized case was obtained in  \cite{SST1}. For more applications about Olsen type inequalities to PDE, one may see \cite{GYT1,GYT2} for details.

 Recently, in
\cite{SST2} the authors obtained an Olsen type inequality for the commutator \ $%
I_{\alpha }^{b}$ \ with a quite elegant method of dyadic decomposition where the definition of the operator $I_{\alpha}^b$ is
 $$
 \begin{array}{ll}
 I_{\alpha}^b(f)(x)=\int_{\Bbb{R}^n}(b(x)-b(y))\frac{f(y)}{|x-y|^{n-\alpha}}dy.
 \end{array}
 \eqno(5.3)
 $$
 with $0<\alpha<n$ and $b\in \text{BMO}(\Bbb{R}^n)$. \\

For the Olsen type inequality of $B\mathcal{I}_\alpha$,  we would like to mention that if we take $v=h$ and $\vec{\omega}=(1,1,\cdots,1)$ in Theorem 4.2, we may obtain
\medskip

 \textbf{Theorem 5.1.} Under the same conditions as in Theorem 4.2, there is

 \textbf{Case 1.} For the  $q>1$, we have
 $$\|h\cdot B\mathcal{I}_\alpha(f,g)\|_{\mathcal{M}^{q_0}_q}\leq C\|h\|_{\mathcal{M}^{r_0}_{r_1}}\|(f,g)\|_{\mathcal{M}^{p_0}_{\vec{P}}}\leq C\|h\|_{\mathcal{M}^{r_0}_{r_1}}\|f\|_{\mathcal{M}_{p_1}^{q_1}}\|g\|_{\mathcal{M}_{p_2}^{q_2}},\eqno(5.4)$$
for all $h\in \mathcal{M}^{r_0}_{r_1}$, $1/q_1+1/q_2=1/p_0$ and $r_1=aq$.

\textbf{Case 2.} For the case $0<q\leq 1$, we have
$$\|h\cdot B\mathcal{I}_\alpha(f,g)\|_{\mathcal{M}^{q_0}_q}\leq C\|h\|_{\mathcal{M}^{r_0}_{q}}\|\vec{f}\|_{\mathcal{M}^{p_0}_{\vec{P}}}\leq C\|h\|_{\mathcal{M}^{r_0}_{q}}\|f\|_{\mathcal{M}_{p_1}^{q_1}}\|g\|_{\mathcal{M}_{p_2}^{q_2}},\eqno(5.5)$$
for all $h\in \mathcal{M}^{r_0}_{q}$ and $1/q_1+1/q_2=1/p_0$.

According to the conditions of Theorem 4.2,  we find that the exponent $r_1=aq$ in (5.4) should satisfy  the condition $r_1\in (q,q\cdot\min\{\frac{r_0}{q_0}, \frac{p_1}{s'},\frac{p_2}{r'}\})\subsetneq (q,r_0)$. Then, comparing (5.2) with  (5.4),   it is natural to ask whether we can get the  following Olsen type inequality for $B\mathcal{I}_\alpha$,
$$\|h\cdot B\mathcal{I}_\alpha(f,g)\|_{\mathcal{M}^{q_0}_q}\leq C\|h\|_{\mathcal{M}^{r_0}_{r_1}}\|f\|_{\mathcal{M}_{p_1}^{q_1}}\|g\|_{\mathcal{M}_{p_2}^{q_2}}$$
with any $r_1\in (q,r_0]$ and $q>1$. In this section, we will give a positive answer to this question. The main result of this section is
\medskip

\textbf{Theorem 5.2.} Suppose that there exist  real  numbers $\alpha,  q_i, p_i \ (i=1,2),  r_0,$\\
 $r_1, s, q_0$ and $q$ satisfying $0<\alpha<n, 1<q_i\leq p_i<\infty, 1<q\leq q_0<\infty,  1<r_1\leq r_0, p_1>r>1, p_2>s>1$ and
$$r_1>q,  1/r_0<\alpha/n<1/q_1+1/q_2<1, 1/s+1/r=1.$$
Furthermore, we assume that
$$1/q_0=1/r_0+1/q_1+1/q_2-\alpha/n$$
and
 $$\frac{q}{q_0}=\frac{p_1}{q_1}=\frac{p_2}{q_2}.\eqno(5.6)$$
Then, there exists a positive constant $C$ independent of $f$ and $g$, such that
$$\| h\cdot B\mathcal{I}_{\alpha}(f,g)\|_{\mathcal{M}^{q_0}_q}\leq  C\| h\|_{\mathcal{M}^{r_0}_{r_1}}\| f\|_{\mathcal{M}^{q_1}_{p_1}}\| g\|_{\mathcal{M}^{q_2}_{p_2}}$$
for any $h\in \mathcal{M}^{r_0}_{r_1}(\Bbb{R}^n)$.

 The method for the proof of Theorem 5.2 is also adapted  to the case $q=\infty$ and $h\equiv 1$. We obtain the Spanne type estimates for $B\mathcal{I}_\alpha$ and it is also a new result with its independent interest.
\medskip

 \textbf{Corollary 5.3. (The Spanne type estimate for $B\mathcal{I}_{\alpha}$)}  Suppose that there exist  real  numbers $\alpha,  q_i, p_i \ (i=1,2), r, s,  q_0 $ and $q$ satisfying $0<\alpha<n, 1<p_i\leq q_i<\infty,  1<q\leq q_0<\infty, p_1>r>1, p_2>s>1$ and
$$  \alpha/n<1/q_1+1/q_2<1, \ \ 1/s+1/r=1.$$
Furthermore, we assume that
$$1/q_0= 1/q_1+1/q_2-\alpha/n$$
and
$$ \frac{q_0}{q}=\frac{q_1}{p_1}=\frac{q_2}{p_2}.$$
Then, there exists a positive constant $C$ independent of $f$ and $g$, such that
$$\|  B\mathcal{I}_{\alpha}(f,g)\|_{\mathcal{M}^{q_0}_q}\leq C\| f\|_{\mathcal{M}^{q_1}_{p_1}}\| g\|_{\mathcal{M}^{q_2}_{p_2}}.$$

\medskip

 \textbf{Remark 5.4.} Here we would like to mention that Theorems 5.2 is  not an  easy consequence of Corollary 5.3 and the H\"{o}lder inequality for functions on the Morrey spaces  (see (2.1) in \cite[p.1377]{ISST2}).  Readers may see \cite{SST3, SST1,YL} for details. In fact, from Corollary 5.3 and the H\"{o}lder inequality for functions on the Morrey spaces, there is
 $$\| h\cdot B\mathcal{I}_{\alpha}(f,g)\|_{\mathcal{M}^{q_0}_q}\leq C\| h\|_{\mathcal{M}^{r_0}_{r_1}}\| f\|_{\mathcal{M}_{p_1}^{q_1}}\| g\|_{\mathcal{M}_{p_2}^{q_2}},\eqno(5.7)$$
 where $$\frac{r_0}{r_1}=\frac{q_0}{q}=\frac{q_1}{p_1}=\frac{q_2}{p_2}\eqno(5.8)$$ and the other conditions are   the same as in Theorem 5.2.
\medskip

\textbf{Remark 5.5.} Comparing (5.6) with (5.8), we find that the restriction of (5.8) is much more stronger than (5.6).
\medskip

\textbf{Remark 5.6.} In \cite{FG},  Fan and Gao \cite[Corollary 2.5]{FG} got an Olsen type inequality for $B\mathcal{I}_{\alpha}$ which is similar to (5.7). If we check \cite[Corollary 2.5]{FG} carefully,  we find that   the exponents $q, q_0, r_1, r_0,  q_1, p_1, q_2, p_2$ in \cite[Corollary 2.5]{FG} also satisfy (5.8). However, our result shows that the condition
 (5.8) is unnecessary as the method used in this paper is quite different and more difficult from \cite{FG}.

\medskip

\textbf{Proof of Theorem 5.2.} Without loss of generality, we may assume that both $f$ and $g$ are non-negative functions.  From Lemma 2.4 and the fact $q\leq q_0$, then  for any cube $Q\subset \Bbb{R}^n$, there is
$$
\begin{array}{ll}
&|Q|^{1/q_0-1/q}\left(\int_{Q}|h(x)B\mathcal{I}_\alpha(f,g)(x)|^qdx\right)^{1/q}\\
&\leq 6^n\sum\limits_{t=1}^{3^n}|Q_t|^{1/q_0-1/q}\left(\int_{Q_t}|h(x)B\mathcal{I}_\alpha(f,g)(x)|^qdx\right)^{1/q},
\end{array}
\eqno(5.9)
$$
where $Q_t\in \mathcal{D}^t$, $Q\subset Q_t$ and $l(Q_t)\leq 6l(Q)$.

 Thus, we only need to estimate $|Q_0|^{1/q_0-1/q}\left(\int_{Q_0}|h(x)B\mathcal{I}_\alpha(f,g)(x)|^qdx\right)^{1/q}$ with $Q_0\in \mathcal{D}^t$.

From (ii) in Section 2, we know that for a fixed $t$ and each $\nu\in \Bbb{Z}$, the set $\{Q\in \mathcal{D}^t: l(Q)=2^{-\nu}\}$ forms a partition of $\Bbb{R}^n$. Moreover, we denote $Q\in \mathcal{D}^t_{\nu}$ with $l(Q)=2^{-\nu}$ and  let $3Q$ be made up of $3^n$ dyadic cubes of equal size and have the same center of $Q$.
Thus,  by the notations as in Section 2, we decompose  $B\mathcal{I}_{ \alpha}$ as follows.

\begin{align*}
&B\mathcal{I}_{\alpha}(f,g)(x)=\int_{\Bbb{R}^n}\frac{f(x-y)g(x+y)}{|y|^{n-\alpha}}dy\\
&=\sum\limits_{\nu\in \Bbb{Z}}\int_{2^{-\nu-1}<|y|\leq 2^{-\nu}}\frac{f(x-y)g(x+y)}{|y|^{n-\alpha}}dy\\
&\leq\sum\limits_{\nu\in \Bbb{Z}}\sum\limits_{Q\in \mathcal{D}^t_\nu}2^{\nu(n-\alpha)}\chi_Q(x)\int_{2^{-\nu-1}<|y|\leq 2^{-\nu}}f(x-y)g(x+y)dy.
\end{align*}
Then, by a geometric observation, we have $B(x,2^{-\nu})\subset 3Q$ if $x\in Q\in \mathcal{D}^t_\nu$. Thus,  using the H\"{o}lder inequality  with $1/r+1/s=1\ \ (r,s>1)$ and a change of variables, there is
\begin{align*}
&\int_{2^{-\nu-1}<|y|\leq 2^{-\nu}}f(x-y)g(x+y)dy\\
&\leq \left(\int_{2^{-\nu-1}<|y|\leq 2^{-\nu}}|f(x-y)|^rdy\right)^{1/r}\left(\int_{2^{-\nu-1}<|y|\leq 2^{-\nu}}|g(x+y)|^sdy\right)^{1/s}\\
&\leq \left(\int_{2^{-\nu-1}<|x-u|\leq 2^{-\nu}}|f(u)|^rdu\right)^{1/r}\left(\int_{2^{-\nu-1}<|x-z|\leq 2^{-\nu}}|g(z)|^sdz\right)^{1/s}\\
&\leq \left(\int_{B(x,2^{-\nu})}|f(u)|^rdu\right)^{1/r}\left(\int_{B(x,2^{-\nu})}|g(z)|^sdz\right)^{1/s}\\
&\leq \left(\int_{3Q}|f(u)|^rdu\right)^{1/r}\left(\int_{3Q}|g(z)|^sdz\right)^{1/s}.
\end{align*}
Then, for any cube fixed cube $Q_0\in \mathcal{D}^t$,  as $Q\in \mathcal{D}^t_\nu$, we  denote

\begin{align*}
I&=h(x)\sum\limits_{\nu\in \Bbb{Z}}\sum\limits_{Q\in \mathcal{D}^t_\nu, Q\supset Q_0}\chi_Q(x)2^{\nu(n-\alpha)} \left(\int_{3Q}|f(u)|^rdu\right)^{1/r}\left(\int_{3Q}|g(z)|^sdz\right)^{1/s}
\end{align*}
and
\begin{align*}
II&=h(x)\sum\limits_{\nu\in \Bbb{Z}}\sum\limits_{Q\in \mathcal{D}^t_\nu, Q\subset Q_0}\chi_Q(x)2^{\nu(n-\alpha)} \left(\int_{3Q}|f(u)|^rdu\right)^{1/r}\left(\int_{3Q}|g(z)|^sdz\right)^{1/s}.
\end{align*}
Thus, it is easy to see
$$h(x)\cdot B\mathcal{I}_\alpha(f,g)(x)\leq I+II.$$
For $I$, let  $Q_k$ be the unique cube containing $Q_0$ and satisfying $|Q_k|=2^{kn}|Q_0|$. Set  $\nu=-\text{log}_2|Q_k|^{\frac{1}{n}}$. Then, we denote
\begin{align*}
E_k&=|Q_0|^{1/q_0-1/q}\left\{ \int_{Q_0}\left|2^{\nu (n-\alpha)}\chi_{Q_k}(x)h(x)\left(\int_{3Q_k}|f(u)|^rdu\right)^{1/r}\right.\right.\\
&\left.\left.\times\left(\int_{3Q_k}|g(z)|^sdz\right)^{1/s}\right|^qdx\right\}^{1/q}.
\end{align*}
Next, we will give the  estimates of $E_k$. By the definition of the Morrey space and the condition $1/r+1/s=1$ with $r,s>1$, we see that
\begin{align*}
&\left(\int_{3Q_k} |f(u)|^rdu\right)^{1/r}\left(\int_{3Q_k} |g(z)|^sdz\right)^{1/s}\\
&\leq \left(\int_{3Q_k}|f(u)|^{p_1}du\right)^{1/p_1}|3Q_k|^{1/r-1/p_1} \left(\int_{3Q_k}|g(z)|^{p_2}dz\right)^{1/p_2}|3Q_k|^{1/s-1/p_2}\\
&\leq \|f\|_{\mathcal{M}^{q_1}_{p_1}}|3Q_k|^{1/r-1/p_1-1/q_1+1/p_1}\|g\|_{\mathcal{M}^{q_2}_{p_2}}|3Q_k|^{1/s-1/p_2+1/p_2-1/q_2}\\
&\leq \|f\|_{\mathcal{M}^{q_1}_{p_1}}\|g\|_{\mathcal{M}^{q_2}_{p_2}}|3Q_k|^{1-1/q_1-1/q_2}.
\end{align*}
Thus, we obtain
\begin{align*}
E_k&\leq \|f\|_{\mathcal{M}_{p_1}^{q_1}}\|g\|_{\mathcal{M}_{p_2}^{q_2}}|3Q_k|^{1-1/q_1-1/q_2}|Q_0|^{1/q_0-1/q}2^{\nu(n-\alpha)}\left(\int_{Q_0}|h(x)|^qdx\right)^{1/q}\\
&\leq \|f\|_{\mathcal{M}_{p_1}^{q_1}}\|g\|_{\mathcal{M}_{p_2}^{q_2}}|3Q_k|^{1-1/q_1-1/q_2}|Q_0|^{1/q_0-1/q+1/q-1/r_1}2^{\nu(n-\alpha)}\left(\int_{Q_0}|h(x)|^{r_1}dx\right)^{1/r_1}\\
&\leq \|h\|_{\mathcal{M}^{r_0}_{r_1}}\|f\|_{\mathcal{M}_{p_1}^{q_1}}\|g\|_{\mathcal{M}_{p_2}^{q_2}}|Q_0|^{1/q_0-1/r_0}|Q_k|^{1-1/q_1-1/q_2}2^{\nu(n-\alpha)}.
\end{align*}
By the facts $2^{\nu(n-\alpha)}=\left(2^{-\text{log}_2|Q_k|^{\frac{1}{n}}}\right)^{n-\alpha}=|Q_k|^{-\frac{1}{n}(n-\alpha)}=|Q_k|^{\frac{\alpha}{n}-1}$ and  $|Q_k|=2^{kn}|Q_0|$, we get
\begin{align*}
E_k&\leq \|h\|_{\mathcal{M}^{r_0}_{r_1}}\|f\|_{\mathcal{M}_{p_1}^{q_1}}\|g\|_{\mathcal{M}_{p_2}^{q_2}}|Q_0|^{1/q_0-1/r_0+\alpha/n-1/q_1-1/q_2}2^{kn(\alpha/n-1/q_1-1/q_2)}\\
&\leq\|h\|_{\mathcal{M}^{r_0}_{r_1}}\|f\|_{\mathcal{M}_{p_1}^{q_1}}\|g\|_{\mathcal{M}_{p_2}^{q_2}}2^{kn(\alpha/n-1/q_1-1/q_2)}.
\end{align*}
Recall that $Q_k$ is the unique cube containing $Q_0$. By the condition that $1/q_1+1/q_2-\alpha/n>0$,  and the definitions of $I$ and $E_{k}$, we obtain
$$
|Q_0|^{1/q_0-1/q}\left(\int_{Q_0}|I|^qdx\right)^{1/q}\leq C\|h\|_{\mathcal{M}^{r_0}_{r_1}}\|f\|_{\mathcal{M}_{p_1}^{q_1}}\|g\|_{\mathcal{M}_{p_2}^{q_2}}.\eqno(5.10)
$$
Next,  we recall some notations from Section 3.  For $r,s>1$ with $1/r+1/s=1$, we set
$$\mathcal{D}_0^t(Q_0)\equiv \{Q\in \mathcal{D}^t(Q_0): m_{3Q}(|f|^r,|g|^s)\leq a\}$$
and
$$\mathcal{D}_{k,j}^t(Q_0)\equiv \{Q\in \mathcal{D}^t(Q_0): Q\subset Q_{k,j}, a^{k}<m_{3Q}(|f|^r,|g|^s)\leq a^{k+1}\},$$
where $a$ is the same as in Section 2 and $\mathcal{D}^t(Q_0)\equiv \{Q\in \mathcal{D}^t: Q\subset Q_0\}$.\\
Thus, we have
$$\mathcal{D}^t(Q_0)=\mathcal{D}_0^t(Q_0)\cup \bigcup\limits_{k,j}\mathcal{D}^t_{k,j}(Q_0).$$
By the duality theory, we may choose a function $\omega\in L^{q'}$, such that
$$\left(\int_{Q_0}|II|^qdx\right)^{1/q}\leq 2\int_{Q_0}|II|\omega(x)dx.\eqno(5.11)$$
Thus, we get
$$
\begin{array}{ll}
&\left(\int_{Q_0}|II|^qdx\right)^{1/q}\\
&=\sum\limits_{Q\in \mathcal{D}_0^t(Q_0)}2^{\nu(n-\alpha)}\int_{Q}h(x)\omega(x)dx\left(\int_{3Q}|f(u)|^rdu\right)^{1/r}\left(\int_{3Q}|g(z)|^sdz\right)^{1/s}\\
&+\sum\limits_{k,j}\sum\limits_{Q\in \mathcal{D}_{k,j}^t(Q_0)}2^{\nu(n-\alpha)}\int_{Q}h(x)\omega(x)dx\left(\int_{3Q}|f(u)|^rdu\right)^{1/r}\left(\int_{3Q}|g(z)|^sdz\right)^{1/s}\\
&:=II_1+II_2.
\end{array}
$$
To estimate $II_2$,   using (2.3), Lemma 2.5,  the definition of $\mathcal{D}_{k,j}(Q_0)$, the geometric property of $\mathcal{D}$ and the fact $0<\frac{\alpha}{n}<1$,  there is
\begin{align*}
II_2&\leq \sum\limits_{k,j}\sum\limits_{Q\in \mathcal{D}_{k,j}^t(Q_0)}2^{\nu(n-\alpha)}\int_{Q}h(x)\omega(x)dx|3Q|m_{3Q}(|f|^r,|g|^s)\\
&\leq C\sum\limits_{k,j}\sum\limits_{Q\in \mathcal{D}_{k,j}^t(Q_0)}|Q|^{\frac{\alpha}{n}}m_{3Q}(|f|^r,|g|^s)\int_Qh(x)\omega(x)dx\\
&\leq C\sum\limits_{k,j}\sum\limits_{Q\in \mathcal{D}_{k,j}^t(Q_0)}|Q|^{\frac{\alpha}{n}}m_{3Q}(|f|^r,|g|^s)\frac{|Q|}{|Q|}\int_Qh(x)\omega(x)dx\\
&\leq C\sum\limits_{k,j}\sum\limits_{Q\in \mathcal{D}_{k,j}^t(Q_0)}|Q|^{\frac{\alpha}{n}}m_{3Q}(|f|^r,|g|^s)\int_QM(h\omega)(x)dx\\
&\leq C\sum\limits_{k,j}|Q_{k,j}|^{\frac{\alpha}{n}}m_{3Q_{k,j}}(|f|^r,|g|^s)\int_{Q_{k,j}}M(h\omega)(x)dx\\
&\leq C\sum\limits_{k,j}|Q_{k,j}|^{\frac{\alpha}{n}}m_{3Q_{k,j}}(|f|^r,|g|^s)m_{Q_{k,j}}[M(h\omega)]|Q_{k,j}|\\
&\leq C\sum\limits_{k,j}|Q_{k,j}|^{\frac{\alpha}{n}}m_{3Q_{k,j}}(|f|^r,|g|^s)m_{Q_{k,j}}[M(h\omega)]|E_{k,j}|.\\
\end{align*}
Thus,  for any $\theta$ satisfying $1<q<\theta<r_1$, we have
\begin{align*}
II_2&\leq C\sum\limits_{k,j}|Q_{k,j}|^{\frac{\alpha}{n}}m_{3Q_{k,j}}(|f|^r,|g|^s)\\
&\times|E_{k,j}|\left(m_{Q_{k,j}}((M^{(\theta')}\omega)^{r_1'})\right)^{1/r_1'}\left(m_{Q_{k,j}}((M^{(\theta)}h)^{r_1})\right)^{1/r_1}\\
&= C\sum\limits_{k,j}|Q_{k,j}|^{\frac{\alpha}{n}-1/r_0}m_{3Q_{k,j}}(|f|^r,|g|^s)|E_{k,j}|\left(m_{Q_{k,j}}((M^{(\theta')}\omega)^{r_1'})\right)^{1/r_1'}\\
&\times|Q_{k,j}|^{1/r_0}\left(m_{Q_{k,j}}((M^{(\theta)}h)^{r_1})\right)^{1/{r_1}}\\
&=\sum\limits_{k,j}|Q_{k,j}|^{\frac{\alpha}{n}-1/r_0}m_{3Q_{k,j}}(|f|^r,|g|^s)|E_{k,j}|\left(m_{Q_{k,j}}((M^{(\theta')}\omega)^{r_1'})\right)^{1/r_1'}\\
\end{align*}
\begin{align*}
&\times\left(|Q_{k,j}|^{\frac{\theta}{r_0}}\left(\frac{1}{|Q_{k,j}|}\int_{Q_{k,j}}M(|h|^\theta)(x)^{r_1/\theta}dx\right)^{\theta/r_1}\right)^{1/\theta}\\
&\leq C\sum\limits_{k,j}|Q_{k,j}|^{\frac{\alpha}{n}-1/r_0}m_{3Q_{k,j}}(|f|^r,|g|^s)|E_{k,j}|\left(m_{Q_{k,j}}((M^{(\theta)}\omega)^{r_1'})\right)^{1/r_1'}\\
&\times|Q_{k,j}|^{1/r_0-1/r_1}\left(\int_{Q_{k,j}}|h(x)|^{r_1}dx\right)^{1/r_1}\\
&\leq C\|h\|_{\mathcal{M}^{r_0}_{r_1}}\sum\limits_{k,j}|Q_{k,j}|^{\frac{\alpha}{n}-1/r_0}m_{3Q_{k,j}}(|f|^r,|g|^s)|E_{k,j}|\left(m_{Q_{k,j}}((M^{(\theta')}\omega)^{r_1'})\right)^{1/r_1'},
\end{align*}
where the definition of $M^{(\theta')}\omega$ can be found in Section 2.

Similarly, for the estimates of $II_1$,  there is
$$II_1\leq C\|h\|_{\mathcal{M}^{r_0}_{r_1}}|Q_{0}|^{\frac{\alpha}{n}-1/r_0}m_{3Q_{0}}(|f|^r,|g|^s)|E_{0}|\left(m_{Q_{0}}((M^{(\theta')}\omega)^{r_1'})\right)^{1/r_1'}.$$
Combing the estimates of $II_1$ and $II_2 $ and  recalling the fact that $\{E_0\}\bigcup\{E_{k,j}\}$ forms a disjoint family of decomposition for $Q_0$, the definition of $Q_{k,j}$ and the fact $\alpha/n>1/r_0$, we get
\begin{align*}
&|Q_0|^{1/q_0-1/q}\int_{Q_0}|II|\omega(x)dx\\
&\leq C|Q_0|^{1/q_0-1/q}\| h\|_{\mathcal{M}^{r_0}_{r_1}}\int_{Q_0}M^{(r_1')}(M^{(\theta')}\omega)(x)M_{\beta_1}(|f|^r)(x)^{1/r}M_{\beta_2}(|g|^s)(x)^{1/s}dx\\
&\leq C|Q_0|^{1/q_0-1/q}\| h\|_{\mathcal{M}^{r_0}_{r_1}}\left(\int_{Q_0}M^{(r_1')}(M^{(\theta')}\omega)(x)^{q'}dx\right)^{1/q'}\\
&\times \left(\int_{Q_0}\left(M_{\beta_1}(|f|^r)(x)^{1/r}M_{\beta_2}(|g|^s)(x)^{1/s}\right)^qdx\right)^{1/q},
\end{align*}
where $M_{\beta_i}$ denotes the fractional maximal function and $\beta_1=\alpha_1r-\frac{nr}{2r_0}>0$, $\beta_2=\alpha_2s-\frac{ns}{2r_0}>0$ with $\alpha_1=\alpha_2=\frac{\alpha}{2}$.\\
As $\frac{q'}{\theta'}>1$ and $\frac{q'}{r_1'}>1$, we can easily get  $\left(\int_{Q_0}M^{(r_1')}(M^{(\theta')}\omega)(x)^{q'}dx\right)^{1/q'}\leq C$ and it remains to give the estimate of
$$|Q_0|^{1/q_0-1/q}\left(\int_{Q_0}\left(M_{\beta_1}(|f|^r)(x)^{1/r}M_{\beta_2}(|g|^s)(x)^{1/s}\right)^qdx\right)^{1/q}.$$
By the H\"{o}lder inequality on Morrey spaces and Lemmas 2.6-2.7, there is
\begin{align*}
|Q_0|^{1/q_0-1/q}&\left(\int_{Q_0}\left(M_{\beta_1}(|f|^r)(x)^{1/r}M_{\beta_2}(|g|^s)(x)^{1/s}\right)^qdx\right)^{1/q}\\
\end{align*}
\begin{align*}
&\leq \|M_{\beta_1}(|f|^r)^{1/r}M_{\beta_2}(|g|^s)^{1/s}\|_{\mathcal{M}^{q_0}_{q}}\\
&\leq \|M_{\beta_1}(|f|^r)^{1/r}\|_{\mathcal{M}^{\mu_1}_{\nu_1}}\|M_{\beta_2}(|g|^s)^{1/s}\|_{\mathcal{M}^{\mu_2}_{\nu_2}}\\
&=\|M_{\beta_1}(|f|^r)\|^{1/r}_{\mathcal{M}^{\frac{\mu_1}{r}}_{\frac{\nu_1}{r}}}\|M_{\beta_2}(|g|^s)\|^{1/s}_{\mathcal{M}^{\frac{\mu_2}{s}}_{\frac{\nu_2}{s}}}\\
&\leq C\||f|^r\|^{1/r}_{M_{\frac{p_1}{r}}^{\frac{q_1}{r}}}\||g|^s\|^{1/s}_{M_{\frac{p_2}{s}}^{\frac{q_2}{s}}}=C\|f\|_{\mathcal{M}^{q_1}_{p_1}}\|g\|_{\mathcal{M}^{q_2}_{p_2}},
\end{align*}
where $\frac{\mu_1}{\nu_1}=\frac{\mu_2}{\nu_2}=\frac{q_0}{q}=\frac{q_1}{p_1}=\frac{q_2}{p_2}$,  $\frac{r}{q_1}-\frac{r}{\mu_1}=\frac{r\alpha_1}{n}-\frac{r}{2r_0}=\frac{\beta_1}{n}$  and  $\frac{s}{q_2}-\frac{s}{\mu_2}=\frac{s\alpha_2}{n}-\frac{s}{2r_0}=\frac{\beta_2}{n}$.
Consequently, we have
$$
|Q_0|^{1/q_0-1/q}\int_{Q_0}|II|\omega(x)dx\leq C\| h\|_{\mathcal{M}^{r_0}_{r_1}}\| f\|_{\mathcal{M}_{p_1}^{q_1}}\| g\|_{\mathcal{M}_{p_2}^{q_2}}.\eqno(5.12)
$$
Thus, combing  (5.9)-(5.12)  , we conclude that

$$
\begin{array}{ll}
\|h\cdot B\mathcal{I}_\alpha(f,g)\|_{\mathcal{M}^{q_0}_q}\leq C\| h\|_{\mathcal{M}^{r_0}_{r_1}}\| f\|_{\mathcal{M}_{p_1}^{q_1}}\| g\|_{\mathcal{M}_{p_2}^{q_2}}.
\end{array}
\eqno(5.13)
$$
Consequently, the proof of Theorem 5.2 has been finished.

\end{document}